\documentclass[12pt,leqno]{article}
\usepackage{amsfonts}
\usepackage{amsmath, amsthm, amsfonts, amssymb, color}
\usepackage{mathrsfs}
\renewcommand{\bar}{\overline}
\renewcommand{\hat}{\widehat}
\renewcommand{\tilde}{\widetilde}

\theoremstyle{definition}

\definecolor{wco}{rgb}{0.5,0.2,0.3}

\numberwithin{equation}{section} \theoremstyle{remark}

\newcommand{\ua}{\uparrow}
\newcommand{\da}{\downarrow}

\title{{\bf Distribution dependent SDEs driven by fractional Brownian motions}
}
\author{
{\bf  Xiliang Fan$^{a),b)}$, Xing Huang$^{c)}$, Yongqiang Suo$^{d)}$, Chenggui Yuan$^{d)}$
 }\\
\footnotesize{$^{a)}$School of Mathematics and Statistics, Anhui Normal University, Wuhu 241002, China}\\
\footnotesize{$^{b)}$Fakult\"{a}t f\"{u}r Mathematik, Universit\"{a}t Bielefeld, 33615 Bielefeld, Germany}\\
\footnotesize{$^{c)}$Center for Applied Mathematics, Tianjin University, Tianjin, 300072, China}\\
\footnotesize{$^{d)}$Department of Mathematics, Swansea University, Bay campus, SA1 8EN, UK}\\
\footnotesize{\sf{fanxiliang0515@163.com},\ \sf{xinghuang@tju.edu.cn},\ \sf{971001@swansea.ac.uk},\ \sf{c.yuan@swansea.ac.uk }}\\
}

\begin{document}
\def\R{\mathbb R}
\def\N{\mathbb N}
\def\E{\mathbb E}
 \def\H{\mathbb H}
\def\Q{\mathbb Q}
\def\H{\mathbb{H}}
\def\P{\mathbb{P}}
\def\S{\mathbb{S}}
\def\Y{\mathbb{Y}}
\def\W{\mathbb{W}}
\def\D{\mathbb{D}}

\def\cH{\mathcal{H}}
\def\cS{\mathcal{S}}

\def\sA{\mathscr{A}}
 \def\sB{\mathscr {B}}
 \def\sC{\mathscr {C}}
  \def\sD{\mathscr {D}}
 \def\sF{\mathscr{F}}
\def\sG{\mathscr{G}}
\def\sL{\mathscr{L}}
\def\sP{\mathscr{P}}
\def\sM{\mathscr{M}}
\def\eq{\equation}
\def\beg{\begin}
\def\ep{\epsilon}
\def\ve{\varepsilon}
\def\vp{\varphi}
\def\vr{\varrho}
\def\om{\omega}
\def\Om{\Omega}
\def\si{\sigma}
\def\ff{\frac}
\def\sq{\sqrt}
\def\kk{\kappa}
\def\de{\delta}
\def\<{\langle} \def\>{\rangle}
\def\Ga{\Gamma}
\def\ga{\gamma}
\def\na{\nabla}
\def\be{\beta}
\def\al{\alpha}
\def\pp{\partial}
 \def\ti{\tilde}
\def\1{\lesssim}
\def\ra{\rightarrow}
\def\da{\downarrow}
\def\upa{\uparrow}
\def\l{\ell}
\def\8{\infty}
\def\3{\triangle}
 \def\DD{\Delta}
\def\m{{\bf m}}
\def\B{\mathbf B}
\def\e{\text{\rm{e}}}
\def\la{\lambda}
\def\th{\theta}

\def\d{\text{\rm{d}}}
   \def\ess{\text{\rm{ess}}}
\def\Ric{\text{\rm{Ric}}} \def \Hess{\text{\rm{Hess}}}
 \def\ua{\underline a}
  \def\Ric{\text{\rm{Ric}}}
\def\cut{\text{\rm{cut}}}      \def\alphaa{\mathbf{r}}     \def\r{r}
\def\gap{\text{\rm{gap}}} \def\prr{\pi_{{\bf m},\varrho}}  \def\r{\mathbf r}

\def\Tilde{\tilde} \def\TILDE{\tilde}\def\II
{\mathbb I}
\def\i{{\rm in}}\def\Sect{{\rm Sect}}

\renewcommand{\bar}{\overline}
\renewcommand{\hat}{\widehat}
\renewcommand{\tilde}{\widetilde}

\allowdisplaybreaks
\maketitle
\begin{abstract}
In this paper we study a class of distribution dependent stochastic differential equations driven by fractional Brownian motions with Hurst parameter $H\in(1/2,1)$.
We prove the well-posedness of this type equations, and then establish a general result on the Bismut formula for the Lions derivative by using Malliavin calculus.
As applications, we provide the Bismut formulas of this kind for both non-degenerate and degenerate cases,
and obtain the estimates of the Lions derivative and the total variation distance between the laws of two solutions.
\end{abstract}


\section{Introduction}

Distribution dependent stochastic differential equations (SDEs), also called McKean-Vlasov or mean-field SDEs,
were initiated by Kac \cite{Kac} in the study of the Boltzman equation and the stochastic toy model.
This type SDEs are equations whose coefficients depend on the distribution of the solution,
and their solutions are often called nonlinear diffusion processes whose transition functions depend on both the current state and distribution of processes.
These distributions dependent SDEs can provide a probabilistic representation to the solutions of a class of nonlinear partial differential equations (PDEs), in which a typical example is propagation of chaos phenomenon (see, e.g., \cite{McKean,Sznitman} and references therein).
Hence, the study of distribution dependent SDEs have received increasing attentions, among which we only mention, for examples,
the works in \cite{CD13,CD15,HMC13} for large scale social interactions within the memory of mean-field games,
\cite{Buckdahn&L&Peng&ainer17a,Crisan&McMurray18a,Li18a} for value functions and related PDEs,
\cite{HRW19,HW,Song,Wang18} for (shift) Harnack type inequalities and gradient estimates, and the references therein.

On the other hand, the Bismut formula, also known as Bismut-Elworthy-Li formula,
was initiated in \cite{Bismut84} by Malliavin calculus and then developed in \cite{Elworthy&Li94} by using martingale method,
which is a very effective tool in the analysis of distributional regularity for various stochastic models.
Afterwards, formula of this kind and (shift) Harnack type inequalities for SDEs and SPDEs were obtained via (backward) coupling argument
(see, e.g., the monograph \cite{Wang13a} and references therein).
Recently, in \cite{RW} the Bismut formula for the Lions derivative ($L$-derivative) was established for distribution dependent SDEs with distribution-free noise by using Malliavin calculus, and then applied to the study of estimates on the $L$-derivative and the total variation distance between distributions of solutions with different initial data.
Here, we mention that in \cite{Banos18} the Bismut formula for initial points was also derived for distribution dependent SDEs,
which can be regarded as a special case of the Bismut formula for the $L$-derivative stated in \cite{RW}.
By introducing the intrinsic and Lions derivatives for probability measures on Banach spaces,
the Bismut formula for the $L$-derivative was given for distribution-path dependent SDEs with distribution-free noise in \cite{BRW}.
Adopting method of heat kernel expansion and the technique of freezing distribution,
the Bismut formula as well as estimate of the $L$-derivative for McKean-Vlasov SDEs with distribution dependent noise was obtained in \cite{HW21}.

In this paper, we are concerned with the following class of distribution dependent SDEs in $\R^d$ on the time interval $[0,T]$:
\beg{align}\label{In1}
\d X_t=b(t,X_t,\sL_{X_t})\d t+\si(t,\sL_{X_t})\d B_t^H,\ \ X_0\in L^p(\Omega\ra\R^d,\sF_0,\P),
\end{align}
where $p\geq 1, \sL_{X_t}$ denotes the law of $X_t, B^H$ is a $d$-dimensional fractional Brownian motion with Hurst parameter $H\in(1/2,1)$,
and the stochastic integral can be regarded as the Wiener integral (see Remark \ref{Rem1} below).
Precise assumptions on the coefficients $b:\Omega\times[0,T]\times\R^d\times\sP_\theta(\R^d)\rightarrow\R^d$ and $\si:[0,T]\times\sP_\theta(\R^d)\rightarrow\R^d\otimes\R^d$ will be specified in later sections,
where $\sP_\th(\R^d)$ is the set of probability measures on $\R^d$ with finite $\th$-th moment.
Let us recall that $B^H=(B^{H,1},\cdots,B^{H,d})$ with Hurst parameter $H\in(0,1)$ is a center Gaussian process with the covariance function $\E(B^{H,i}_tB^{H,j}_s)=R_H(t,s)\delta_{i,j}$, where
\beg{align*}
 R_H(t,s)=\frac{1}{2}\left(t^{2H}+s^{2H}-|t-s|^{2H}\right),\ \  t,s\in[0,T].
\end{align*}
This implies that the relation $\E(|B_t^{H,i}-B_s^{H,i}|^q)=C_q|t-s|^{qH}$ holds for every $q\geq 1$ and $i=1,\cdots,d$.
Consequently, $B^H$ is $(H-\ve)$-order H\"{o}lder continuous a.s. for any $\ve\in(0,H)$ and is an $H$-self similar process.
This, together with the fact that $B^{1/2}$ is a standard Brownian motion, converts fractional Brownian motion into a natural generalization of Brownian motion and leads to many applications in modelling physical phenomena and finance behaviours.

The aim of this paper is to establish a Bismut type $L$-derivative formula for the equation \eqref{In1}.
We first prove the well-posedness of \eqref{In1} under $\W_\th$-Lipschitz conditions with respect to the measure variable,
in which the solution $X$ belongs to the space $X\in\cS^p([0,T])$  that generalizes and improves the corresponding one in the existing literature (see Theorem \ref{(EX)The1} and Remark \ref{WP-rem} below).
Then, with the help of partial derivatives with respect to the initial value and the Malliavin derivative of solutions to \eqref{In1} obtained under stronger assumptions than that for well-posedness, we are able to establish a general result concerning the Bismut type $L$-derivative formula for \eqref{In1} (see Theorem \ref{(Ge)Th1} below).
As applications, we provide the Bismut type $L$-derivative formulas for both non-degenerate and degenerate situations.
In addition, to illustrate the power of the Bismut type $L$-derivative formulas, in the case of non-degenerate we obtain the estimates of the $L$-derivative and the total variation distance for the difference between the laws of the solutions $\sL_{X_T^\mu}$ and $\sL_{X_T^\nu}$ with different initial distributions  $\mu$ and $\nu$ (see Remark \ref{Non-Re2} (ii) below).

The plan of the paper is as follows:
Section 2 is devoted to recalling some useful facts on fractional calculus, fractional Brownian motion and the $L$-derivative.
In Section 3, we prove the existence and uniqueness of a solution to distribution dependent SDE driven by fractional Brownian motion.
In Section 4, we state and prove our main results concerning Bismut type formula for the $L$-derivative of distribution dependent SDE driven by fractional Brownian motion, which are then applied to both non-degenerate and degenerate cases.

\section{Preliminaries}

This section is devoted to giving some basic elements of fractional calculus involving fractional integral and derivative,
Wiener space associated to fractional Brownian motion and the Lions derivative.

\subsection{Fractional integral and derivative}

Let $a,b\in\R$ with $a<b$.
For $f\in L^1([a,b],\R)$ and $\alpha>0$, the left-sided (respectively right-sided) fractional Riemann-Liouville integral of $f$ of order $\alpha$
on $[a,b]$ is defined as
\beg{align*}
&I_{a+}^\alpha f(x)=\frac{1}{\Gamma(\alpha)}\int_a^x\frac{f(y)}{(x-y)^{1-\alpha}}\d y\\
&\qquad\left(\mbox{respectively}\ \ I_{b-}^\alpha f(x)=\frac{(-1)^{-\alpha}}{\Gamma(\alpha)}\int_x^b\frac{f(y)}{(y-x)^{1-\alpha}}\d y\right),
\end{align*}
where $x\in(a,b)$ a.e., $(-1)^{-\alpha}=\e^{-i\alpha\pi}$ and $\Gamma$ denotes the Gamma function.
In particular, when $\alpha=n\in\N$, they are consistent with the usual $n$-order iterated integrals.

Fractional differentiation may be given as an inverse operation.
Let $\alpha\in(0,1)$ and $p\geq1$.
If $f\in I_{a+}^\alpha(L^p([a,b],\R))$ (respectively $I_{b-}^\alpha(L^p([a,b],\R)))$, then the function $g$ satisfying $f=I_{a+}^\alpha g$ (respectively $f=I_{b-}^\alpha g$) is unique in $L^p([a,b],\R)$ and it coincides with the left-sided (respectively right-sided) Riemann-Liouville derivative
of $f$ of order $\alpha$ shown by
\beg{align*}
&D_{a+}^\alpha f(x)=\frac{1}{\Gamma(1-\alpha)}\frac{\d}{\d x}\int_a^x\frac{f(y)}{(x-y)^\alpha}\d y\\
&\qquad\left(\mbox{respectively}\ D_{b-}^\alpha f(x)=\frac{(-1)^{1+\alpha}}{\Gamma(1-\alpha)}\frac{\d}{\d x}\int_x^b\frac{f(y)}{(y-x)^\alpha}\d y\right).
\end{align*}
The corresponding Weyl representation is of the form
\beg{align}\label{FrDe}
&D_{a+}^\alpha f(x)=\frac{1}{\Gamma(1-\alpha)}\left(\frac{f(x)}{(x-a)^\alpha}+\alpha\int_a^x\frac{f(x)-f(y)}{(x-y)^{\alpha+1}}\d y\right)\\
&\qquad\left(\mbox{respectively}\ \ D_{b-}^\alpha f(x)=\frac{(-1)^\alpha}{\Gamma(1-\alpha)}\left(\frac{f(x)}{(b-x)^\alpha}+\alpha\int_x^b\frac{f(x)-f(y)}{(y-x)^{\alpha+1}}\d y\right)\right),\nonumber
\end{align}
where the convergence of the integrals at the singularity $y=x$ holds pointwise for almost all $x$ if $p=1$ and in the $L^p$ sense if $p>1$.
For further details, we refer the reader to \cite{SK}.

\subsection{Wiener space associated to fractional Brownian motion}

For fixed $H\in(1/2,1)$, let $(\Omega,\sF,\P)$ be the canonical probability space associated with fractional Brownian motion with Hurst parameter $H$.
More precisely,
$\Omega$ is the Banach space $C_0([0,T],\R^d)$ of continuous functions vanishing at $0$ equipped with the supremum norm,
$\sF$ is the Borel $\si$-algebra and $\P$ is the unique probability measure on $\Omega$ such that the canonical process $\{B^H_t; t\in[0,T]\}$ is a $d$-dimensional fractional Brownian motion with Hurst parameter $H$.
We assume that there is a sufficiently rich sub-$\si$-algebra $\sF_0\subset\sF$ independent of $B^H$ such that
for any $\mu\in\sP_2(\R^d)$ there exists a random variable $X\in L^2(\Omega\ra\R^d,\sF_0,\P)$ with distribution $\mu$.
Let $\{\sF_t\}_{t\in[0,T]}$ be the filtration generated by $B^H$, completed and augmented by $\sF_0$.

Let $\mathscr{E}$ be the set of step functions on $[0,T]$ and $\mathcal {H}$ the Hilbert space defined as the closure of
$\mathscr{E}$ with respect to the scalar product
\beg{align*}
\left\langle (\mathbb I_{[0,t_1]},\cdot\cdot\cdot,\mathbb I_{[0,t_d]}),(\mathbb I_{[0,s_1]},\cdot\cdot\cdot,\mathbb I_{[0,s_d]})\right\rangle_\cH=\sum\limits_{i=1}^dR_H(t_i,s_i).
\end{align*}
The mapping $(\mathbb I_{[0,t_1]},\cdot\cdot\cdot,\mathbb I_{[0,t_d]})\mapsto\sum_{i=1}^dB_{t_i}^{H,i}$ can be extended to an isometry between $\cH$ (also called the reproducing kernel Hilbert space) and the Gaussian space $\mathcal {H}_1$ associated with $B^H$.
Denote this isometry by $\psi\mapsto B^H(\psi)$.
On the other hand, it follows from \cite{Decreusefond&Ustunel98a} that $R_H(t,s)$ has the following integral representation
\beg{align*}
 R_H(t,s)=\int_0^{t\wedge s}K_H(t,r)K_H(s,r)\d r,
\end{align*}
where $K_H$ is a square integrable kernel given by
\beg{align*}
K_H(t,s)=\Gamma\left(H+\frac{1}{2}\right)^{-1}(t-s)^{H-\frac{1}{2}}F\left(H-\frac{1}{2},\frac{1}{2}-H,H+\frac{1}{2},1-\frac{t}{s}\right),
\end{align*}
in which $F(\cdot,\cdot,\cdot,\cdot)$ is the Gauss hypergeometric function (for details see \cite{Decreusefond&Ustunel98a} or \cite{Nikiforov&Uvarov88}).

Now, define the linear operator $K_H^*:\mathscr{E}\rightarrow L^2([0,T],\R^d)$ as follows
\beg{align*}
(K_H^*\psi)(s)=K_H(T,s)\psi(s)+\int_s^T(\psi(r)-\psi(s))\frac{\partial K_H}{\partial r}(r,s)\d r.
\end{align*}
According to \cite{Alos&Mazet&Nualart01a},
the relation $\langle K_H^*\psi,K_H^*\phi\rangle_{L^2([0,T],\R^d)}=\langle\psi,\phi\rangle_\mathcal {H}$ holds for all $\psi,\phi\in\mathscr{E}$,
and then by the bounded linear transform theorem, $K_H^*$ can be extended to an isometry between $\mathcal{H}$ and $L^2([0,T],\R^d)$.
Consequently, by \cite{Alos&Mazet&Nualart01a} again, the process $\{W_t=B^H((K_H^*)^{-1}{\mathbb I}_{[0,t]}),t\in[0,T]\}$ is a Wiener process,
and $B^H$ has the following Volterra-type representation
\beg{align}\label{IRFor}
B_t^H=\int_0^tK_H(t,s)\d W_s, \ \ t\in[0,T].
\end{align}
Besides, we define the operator $K_H: L^2([0,T],\mathbb{R}^d)\rightarrow I_{0+}^{H+1/2 }(L^2([0,T],\mathbb{R}^d))$ by
\beg{align*}
 (K_H f)(t)=\int_0^tK_H(t,s)f(s)\d s.
\end{align*}
Due to \cite{Decreusefond&Ustunel98a}, we know that it is an isomorphism and for each $f\in L^2([0,T],\mathbb{R}^d)$,
\beg{align*}
 (K_H f)(s)=I_{0+}^{1}s^{H-1/2}I_{0+}^{H-1/2}s^{1/2-H}f.
\end{align*}
Then for every $h\in I_{0+}^{H+1/2}(L^2([0,T],\R^d))$, the inverse operator $K_H^{-1}$ is of the form
\beg{align}\label{InOp}
(K_H^{-1}h)(s)=s^{H-1/2}D_{0+}^{H-1/2}s^{1/2-H}h'.
\end{align}
We remark that the injection $R_H=K_H\circ K_H^*:\mathcal{H}\rightarrow\Omega$ embeds $\mathcal{H}$ densely into $\Omega$ and
for every $\psi\in\Omega^*\subset\mathcal{H}$ there holds
$\E\e^{i\<B^H,\psi\>}=\exp(-\ff 1 2\|\psi\|^2_\mathcal{H})$.
Consequently, $(\Omega,\mathcal{H},\P)$ is an abstract Wiener space in the sense of Gross.

Finally, we give a brief account on the Malliavin calculus for fractional Brownian motion.
Denote $\mathcal {S}$ by the set of smooth and cylindrical random variables of the form
\beg{align*}
F=f(B^H(\phi_1),\cdot\cdot\cdot,B^H(\phi_n)),
\end{align*}
where $n\geq 1, f\in C_b^\infty(\mathbb{R}^n)$, which is the collection of $f$ and all its partial derivatives are bounded,
$\phi_i\in\mathcal{H}, 1\leq i\leq n$.
The Malliavin derivative of $F$, denoted by $\D F$, is defined as the $\mathcal {H}$-valued random variable
\beg{align*}
\D F=\sum_{i=1}^n\frac{\partial f}{\partial x_i}(B^H(\phi_1),\cdot\cdot\cdot,B^H(\phi_n))\phi_i.
\end{align*}
For any $p\geq 1$, we define the Sobolev space $\mathbb{D}^{1,p}$ as the completion of $\mathcal {S}$ with respect to the norm
\beg{align*}
\|F\|_{1,p}^p=\mathbb{E}|F|^p+\mathbb{E}\|\D F\|^p_{\mathcal {H}}.
\end{align*}
Meanwhile, we will denote by $\delta$ and $\mathrm{Dom}\delta$ the dual operator of $\mathbb{D}$ and its domain, respectively.
Let us finish this part by giving a transfer principle that connects the derivative and divergence operators of both processes $B^H$ and $W$ that are needed later on.
\beg{prp}\label{Tra-Prp1}\cite[Proposition 5.2.1]{ND}
For any $F\in\D_W^{1,2}=\D^{1,2}$,
\beg{align*}
K_H^*\D F=\D^WF,
\end{align*}
where $\D^W$ denotes the derivative operator with respect to the underlying Wiener process $W$ appearing in \eqref{IRFor}, and $\D_W^{1,2}$ the corresponding Sobolev space.
\end{prp}

\beg{prp}\label{Tra-Prp2}\cite[Proposition 5.2.2]{ND}
$\mathrm{Dom}\delta=(K_H^*)^{-1}(\mathrm{Dom}\delta_W)$,
and for any $\mathcal{H}$-valued random variable $u$ in $\mathrm{Dom}\delta$ we have $\delta(u)=\delta_W(K_H^*u)$,
where $\delta_W$ denotes the divergence operator with respect to the underlying Wiener process $W$ appearing in \eqref{IRFor}.
\end{prp}

\beg{rem}\label{Tra-Rem1}
The above proposition, together with \cite[Proposition 1.3.11]{ND}, yields that
if $K_H^*u\in L_a^2([0,T]\times\Omega,\R^d)$ (the closed subspace of $L^2([0,T]\times\Omega,\R^d)$ formed by the adapted processes),
then $u\in\mathrm{Dom}\delta$.
\end{rem}

\subsection{The Lions derivative}

For $\th\in[1,\infty)$, let $\sP_\th(\R^d)$ be the space of probability measures on $\R^d$ with finite $\th$-th moment.
We define the $L^\th$-Wasserstein distance on $\sP_\th(\R^d)$ by
$$\mathbb{W}_\th(\mu,\nu):=\inf_{\pi\in\sC(\mu,\nu)}\left(\int_{\R^d\times\R^d}|x-y|^\th\pi(\d x, \d y)\right)^\ff 1 \th,\ \ \mu,\nu\in\sP_\th(\R^d),$$
where $\sC(\mu,\nu)$ is the set of probability measures on $\R^d\times\R^d$ with marginals $\mu$ and $\nu$.
It is well-known that $(\sP_\th(\R^d),\mathbb{W}_\th)$ is a Polish space.
Throughout this paper, we use $|\cdot|$ and $\<\cdot,\cdot\>$ for the Euclidean norm and inner product, respectively,
and for a matrix, we denote by $\|\cdot\|$ either the operator norm or the Euclidean norm.
$\|\cdot\|_{L^2_\mu}$ denotes for the norm of the space $ L^2(\R^d\ra\R^d,\mu)$ and for a random variable $X$, $\sL_X$ denotes its distribution.

\beg{defn}
Let $f:\sP_2(\R^d)\ra\R$ and $g:\R^d\times\sP_2(\R^d)\ra\R$.
\beg{enumerate}
\item[(1)]  $f$ is called $L$-differentiable at $\mu\in\sP_2(\R^d)$, if the functional
\beg{align*}
L^2(\R^d\ra\R^d,\mu)\ni\phi\mapsto f(\mu\circ(\mathrm{Id}+\phi)^{-1}))
\end{align*}
is Fr\'{e}chet differentiable at $0\in L^2(\R^d\ra\R^d,\mu)$, i.e., there exists a unique $\gamma\in L^2(\R^d\ra\R^d,\mu)$ such that
\beg{align*}
\lim_{\|\phi\|_{L^2_\mu}\ra0}\ff{f(\mu\circ(\mathrm{Id}+\phi)^{-1})-f(\mu)-\mu(\<\gamma,\phi\>)}{\|\phi\|_{L^2_\mu}}=0.
\end{align*}
In this case, $\gamma$ is called the $L$-derivative of $f$ at $\mu$ and denoted by $D^Lf(\mu)$.

\item[(2)] $f$ is called $L$-differentiable on $\sP_2(\R^d)$, if the $L$-derivative $D^Lf(\mu)$ exists for all $\mu\in\sP_2(\R^d)$.
If, moreover, for every $\mu\in\sP_2(\R^d)$ there exists a $\mu$-version $D^Lf(\mu)(\cdot)$ such that $D^Lf(\mu)(x)$ is jointly continuous in $(\mu,x)\in\sP_2(\R^d)\times\R^d$, we denote $f\in C^{(1,0)}(\sP_2(\R^d))$.

\item[(3)] $g$ is called differentiable on $\R^d\times\sP_2(\R^d)$, if for any $(x,\mu)\in\R^d\times\sP_2(\R^d)$,
$g(\cdot,\mu)$ is differentiable and $g(x, \cdot)$ is $L$-differentiable.
If, moreover, $\na g(\cdot,\mu)(x)$ and $D^Lg(x,\cdot)(\mu)(y)$ are jointly continuous in $(x,y,\mu)\in\R^d\times\R^d\times\sP_2(\R^d)$,
we denote $g\in C^{1,(1,0)}(\R^d\times\sP_2(\R^d)$.

\end{enumerate}
\end{defn}
As mentioned in \cite[Section 2]{RW21}, the above definition of $L$-derivative coincides with the Wasserstein derivative,
which was introduced by P.-L. Lions using probability spaces (\cite{Cardaliaguet13}).
Besides, it is easy to see that if $f$ is  $L$-differentiable at $\mu$, then for any $\phi\in L^2(\R^d\ra\R^d,\mu)$ there holds
\beg{align*}
D^L_\phi f(\mu):=\lim_{\ve\da0}\ff{f(\mu\circ(\mathrm{Id}+\ve\phi)^{-1})-f(\mu)}\ve=\mu(\<D^Lf(\mu),\phi\>),
\end{align*}
in which $D^L_\phi f(\mu)$ is called the directional $L$-derivative of $f$ along $\phi$ initiated in \cite{AKR96}.
When $D^L_\cdot f(\mu):L^2(\R^d\ra\R^d,\mu)\ra\R$ is a bounded linear functional,
$\phi\mapsto f(\mu\circ(\mathrm{Id}+\phi)^{-1}))$ is G\^{a}teaux differentiable at $0$.
In this case, we say that $f$ is weakly $L$-differentiable at $\mu$ (also called intrinsically differentiable at $\mu$, see \cite[Definition 2.1]{BRW} or \cite[Definition 1.1]{HW21}).
Moreover, we have
\beg{align*}
\|D^Lf(\mu)\|:=\|D^Lf(\mu)(\cdot)\|_{L^2_\mu}=\sup_{\|\phi\|_{L^2_\mu}\leq 1}|D^L_\phi f(\mu)|.
\end{align*}
For a vector-valued function $f=(f_i)$ or a matrix-valued function $f=(f_{ij})$ with $L$-differentiable components, we simply write
\beg{align*}
D^Lf(\mu)=(D^Lf_i(\mu))  \ \  \mathrm{or}\ \ D^Lf(\mu)=(D^Lf_{ij}(\mu)).
\end{align*}
Besides, we have the following useful formula for the $L$-derivative, which is due to \cite[Theorem 6.5]{Cardaliaguet13} and \cite[Proposition 3.1]{RW}.
\beg{lem}\label{FoLD}
Let $(\Omega,\sF,\P)$ be an atomless probability space and $(X,Y)\in L^2(\Omega\ra R^d,\P)$.
If $f\in C^{1,0}(\sP_2(\R^d))$, then
\beg{align*}
\lim_{\ve\da0}\ff {f(\sL_{X+\ve Y})-f(\sL_X)} \ve=\E\<D^Lf(\sL_X)(X),Y\>.
\end{align*}
\end{lem}

\section{Well-posedness of distribution dependent SDE by fractional noise}

In this section, we consider the following distribution dependent SDE driven by fractional Brownian motion:
\beg{align}\label{XX}
\d X_t=b(t,X_t,\sL_{X_t})\d t+\si(t,\sL_{X_t})\d B_t^H,\ \ X_0=\xi,
\end{align}
where the coefficients $b:\Omega\times[0,T]\times\R^d\times\sP_\theta(\R^d)\rightarrow\R^d, \si:[0,T]\times\sP_\theta(\R^d)\rightarrow\R^d\otimes\R^d$
and $\xi\in L^p(\Omega\ra\R^d,\sF_0,\P)$ with $p\geq \theta(\geq1)$.
Now, we introduce the hypotheses under which we will be able to prove the existence and uniqueness of solutions to \eqref{XX}.
\beg{enumerate}
\item[(H)] There exists a non-decreasing function $K(t)$ such that for any $t\in[0,T], x,y\in\R^d, \mu,\nu\in\sP_\theta(\R^d)$,
\beg{align*}
|b(t,x,\mu)-b(t,y,\nu)|\le K(t)(|x-y|+\W_\theta(\mu,\nu)),\ \ \|\si(t,\mu)-\si(t,\nu)\|\le K(t)\W_\theta(\mu,\nu),
\end{align*}
and
\beg{align*}
|b(t,0,\de_0)|+\|\si(t,\de_0)\|\le K(t).
\end{align*}
\end{enumerate}
For any $p\geq1$, let $\cS^p([0,T])$ be the space of  $\R^d$-valued, continuous $(\sF_t)_{t\in[0,T]}$-adapted processes $\psi$ on $[0,T]$ satisfying
$$\|\psi\|_{\cS^p}:=\bigg(\E\sup_{t\in[0,T]}|\psi_t|^p\bigg)^{1/p}<\8,$$
and let the letter $C$ with or without indices denote generic constants, whose values may change from line to line.

\beg{defn}\label{Def1}
A stochastic process $X=(X_t)_{0\leq t \leq T}$ on $\R^d$ is called a solution of \eqref{XX},
if $X\in\cS^p([0,T])$ and $\P$-a.s.,
$$X_t=\xi+\int_0^tb(s,X_s,\sL_{X_s})\d s+\int_0^t\si(s,\sL_{X_s})\d B_s^H,\ \ t\in[0,T].$$
\end{defn}

\beg{rem}\label{Rem1}
Observe that $\si(\cdot,\sL_{X_\cdot})$ is a deterministic function,
then $\int_0^t\si(s,\sL_{X_s})\d B_s^H$ is regarded as a Wiener integral with respect to fractional Brownian motion.
\end{rem}

\beg{thm}\label{(EX)The1}
Suppose that (H) holds and $\xi\in L^p(\Omega\ra\R^d,\sF_0,\P)$ with $p\geq \theta$ and $p>1/H$.
Then the equation \eqref{XX} has a unique solution $X\in\cS^p([0,T])$.
\end{thm}

Before proving the theorem, we first present the following Hardy-Littlewood inequality (see, e.g., \cite[Theroem 1, Page 119]{Stein70}).
\beg{lem}\label{HLI}
Let $1<\tilde{p}<\tilde{q}<\infty$ and $\frac{1}{\tilde{q}}=\frac{1}{\tilde{p}}-\alpha$.
If $f:\R_+\rightarrow\R$ belongs to $L^{\tilde{p}}(0,\infty)$, then $I^\alpha_{0+}f(x)$ converges absolutely for almost every $x$, and moreover
$$\|I^\alpha_{0+}f\|_{L^{\tilde{q}}(0,\infty)}\leq C_{\tilde{p},\tilde{q}}\|f\|_{L^{\tilde{p}}(0,\infty)}$$
holds for some positive constant $C_{\tilde{p},\tilde{q}}$.
\end{lem}

\emph{Proof of Theorem \ref{(EX)The1}.}
Define recursively $(X^n)_{n\geq 1}$ as follows:
$X^0_t=\xi,\ t\in[0,T]$ and for each $n\geq 1$,
\beg{align*}
X^n_t=\xi+\int_0^tb(s,X^{n-1}_s,\sL_{X^{n-1}_s})\d s+\int_0^t\si(s,\sL_{X^{n-1}_s})\d B_s^H,\ \ t\in[0,T].
\end{align*}
The rest of the proof will be divided into three steps.

\textit{Step 1. Claim: For any $p\geq \theta$ and $p>1/H$, if $\E\Big(\sup_{t\in[0,T]}|X^n_t|^p\Big)<\infty$, then there holds $\E\Big(\sup_{t\in[0,T]}|X^{n+1}_t|^p\Big)<\infty$}.
Owing to the H\"{o}lder inequality  and  (H), we have for any $p\geq\theta$,
\beg{align}\label{Pf1-The1}
\E\bigg(\sup\limits_{t\in[0,T]}|X^{n+1}_t|^p\bigg)&\leq 3^{p-1}\E|\xi|^p+3^{p-1}\E\left(\sup\limits_{t\in[0,T]}\left|\int_0^tb(s,X^n_s,\sL_{X^n_s})\d s\right|^p\right)\cr
&\ \ \ +3^{p-1}\E\left(\sup\limits_{t\in[0,T]}\left|\int_0^t\si(s,\sL_{X^n_s})\d B_s^H\right|^p\right)\cr
&\leq 3^{p-1}\E|\xi|^p+(3T)^{p-1}\E\int_0^TK^p(s)(1+|X_s^n|+\W_\theta(\sL_{X^n_s},\delta_0))^p\d s\cr
&\ \ \ +3^{p-1}\E\left(\sup\limits_{t\in[0,T]}\left|\int_0^t\si(s,\sL_{X^n_s})\d B_s^H\right|^p\right)\cr
&\leq 3^{p-1}\E|\xi|^p+3^{2(p-1)}(T K(T))^p\left(1+2\E\sup_{t\in[0,T]}|X^n_t|^p\right)\cr
&\ \ \ +3^{p-1}\E\left(\sup\limits_{t\in[0,T]}\left|\int_0^t\si(s,\sL_{X^n_s})\d B_s^H\right|^p\right).
\end{align}

Next, we shall provide an estimate for the last term of the right-hand side of \eqref{Pf1-The1},
whose argument is partially borrowed from \cite[Theorem 4]{AN03}.

We take $\la$ satisfying $1-H<\la<1-1/p$ because $pH>1$.
Using the fact that $\int_s^t(t-r)^{-\la}(r-s)^{\la-1}\d r=C(\la)$, the stochastic Fubini theorem and the H\"{o}lder inequality, we get
\beg{align}\label{Pf2-The1}
&\E\left(\sup\limits_{t\in[0,T]}\left|\int_0^t\si(s,\sL_{X^n_s})\d B_s^H\right|^p\right)\cr
&=C(\la)^{-p}\E\left(\sup\limits_{t\in[0,T]}\left|\int_0^t\left(\int_s^t(t-r)^{-\la}(r-s)^{\la-1}\d r\right)\si(s,\sL_{X^n_s})\d B_s^H\right|^p\right)\cr
&=C(\la)^{-p}\E\left(\sup\limits_{t\in[0,T]}\left|\int_0^t (t-r)^{-\la}\left(\int_0^r(r-s)^{\la-1}\si(s,\sL_{X^n_s})\d B_s^H\right)\d r\right|^p\right)\cr
&\leq\ff {C(\la)^{-p}}{(p-1-\la p)^{p-1}}\E\left(\sup\limits_{t\in[0,T]}t^{p-1-\la p}
\int_0^t \left|\int_0^r(r-s)^{\la-1}\si(s,\sL_{X^n_s})\d B_s^H\right|^p\d r\right)\cr
&\leq\ff {C(\la)^{-p}}{(p-1-\la p)^{p-1}}T^{p-1-\la p}\int_0^T\E\left|\int_0^r(r-s)^{\la-1}\si(s,\sL_{X^n_s})\d B_s^H\right|^p\d r,
\end{align}
where we use the condition $\la<1-1/p$ in the first inequality.\\
Notice that for each $r\in[0,T], \int_0^r(r-s)^{\la-1}\si(s,\sL_{X^n_s})\d B_s^H$ is a centered Gaussian random variable.
Then by the Kahane-Khintchine formula, we obtain that there exists a constant $C_p>0$ such that
\beg{align}\label{Pf3-The1}
&\E\left|\int_0^r(r-s)^{\la-1}\si(s,\sL_{X^n_s})\d B_s^H\right|^p\cr
&\leq C_p\left(\E\left|\int_0^r(r-s)^{\la-1}\si(s,\sL_{X^n_s})\d B_s^H\right|^2\right)^{\ff p 2}\cr
&\leq C_p\left(\int_0^r\int_0^r(r-u)^{\la-1}\|\si(u,\sL_{X^n_u})\|(r-v)^{\la-1}\|\si(v,\sL_{X^n_v})\||u-v|^{2H-2}\d u\d v\right)^{\ff p 2}\cr
&\leq C_{p,H}\left(\int_0^r(r-s)^{\ff {\la-1} H} \|\si(s,\sL_{X^n_s})\|^{\ff 1 H}\d s\right)^{pH},
\end{align}
where the last inequality is due to the argument in \cite[Theroem 1.1, Page 201]{JMV01}.\\
Substituting \eqref{Pf3-The1} into \eqref{Pf2-The1} and using the condition $1-H<\la$ and Lemma \ref{HLI} with $\tilde{q}=pH$ and $\al=1-\ff{1-\la}H$
(imply $\tilde{p}=\ff{pH}{p(\la+H-1)+1}$), we have
\beg{align}\label{MI}
&\E\left(\sup\limits_{t\in[0,T]}\left|\int_0^t\si(s,\sL_{X^n_s})\d B_s^H\right|^p\right)\cr
&\leq C_{\la,p,H}T^{p-1-\la p}\int_0^T\left(\int_0^r(r-s)^{\ff {\la-1} H} \|\si(s,\sL_{X^n_s})\|^{\ff 1 H}\d s\right)^{pH}\d r\cr
&\leq C_{\la,p,H}T^{p-1-\la p}\left(\int_0^T \|\si(r,\sL_{X^n_r})\|^{\ff p {p(\la+H-1)+1}}\d r\right)^{p(\la+H-1)+1}\cr
&\leq C_{\la,p,H}T^{pH-1}\int_0^T\|\si(s,\sL_{X^n_s})\|^p\d s,
\end{align}
where we use the H\"{o}lder inequality in the last inequality, and remark that $C_{\la,p,H}$ above may depend only on $p$ and $H$ by choosing proper $\la$.\\
Observe that, by (H) and $p\geq\theta$ we have
\beg{align*}
\int_0^T\|\si(s,\sL_{X^n_s})\|^p\d s&\leq\int_0^TK^p(s)\left(1+\W_\theta(\sL_{X^n_s},\delta_0)\right)^p\d s\cr
&\leq 2^{p-1}K^p(T)T\left(1+\E\bigg(\sup_{t\in[0,T]}|X^n_t|^p\bigg)\right).
\end{align*}
Then, plugging this into \eqref{MI} yields
\beg{align*}
\E\left(\sup\limits_{t\in[0,T]}\left|\int_0^t\si(s,\sL_{X^n_s})\d B_s^H\right|^p\right)
\leq C_{p,H}K^p(T)T^{pH}\left(1+\E\bigg(\sup_{t\in[0,T]}|X^n_t|^p\bigg)\right).
\end{align*}
Combining this with \eqref{Pf1-The1} and the assumption that $\E(\sup_{t\in[0,T]}|X^n_t|^p)<\infty$ yields the desired claim.

\textit{Step 2. Existence.}
To this end, we shall prove the convergence of $X^n$ in $\cS^p([0,T])$ with any $p\geq\theta$.
For any $t\in[0,T]$, we get
\beg{align}\label{Pf4-The1}
\E\bigg(\sup\limits_{s\in[0,t]}|X^n_s-X^{n-1}_s|^p\bigg)&\leq 2^{p-1}\E\left(\sup\limits_{s\in[0,t]}\left|\int_0^s\left(b(r,X^{n-1}_r,\sL_{X^{n-1}_r})-b(r,X^{n-2}_r,\sL_{X^{n-2}_r})\right)\d r\right|^p\right)\nonumber\\
&\ \ \ +2^{p-1}\E\left(\sup\limits_{s\in[0,t]}\left|\int_0^s\left(\si(r,\sL_{X^{n-1}_r})-\si(r,\sL_{X^{n-2}_r})\right)\d B_r^H\right|^p\right)\nonumber\\
&=:2^{p-1}I_1(t)+2^{p-1}I_2(t).
\end{align}
For the term $I_1(t)$, from (H) and $p\geq\theta$ we obtain
\beg{align}\label{Pf5-The1}
I_1(t)&\leq t^{p-1}\E\left(\sup\limits_{s\in[0,t]}\int_0^s\left|b(r,X^{n-1}_r,\sL_{X^{n-1}_r})-b(r,X^{n-2}_r,\sL_{X^{n-2}_r})\right|^p\d r\right)\cr
&\leq t^{p-1}\E\left(\sup\limits_{s\in[0,t]}\int_0^s\left[K(r)\left(|X^{n-1}_r-X^{n-2}_r|+\W_\theta(\sL_{X^{n-1}_r},\sL_{X^{n-2}_r})\right)\right]^p\d r\right)\cr
&\leq 2^{p-1}K^p(t)t^{p-1}\E\left(\sup\limits_{s\in[0,t]}\int_0^s\left(|X^{n-1}_r-X^{n-2}_r|^p+\E|X^{n-1}_r-X^{n-2}_r|^p\right)\d r\right)\cr
&\leq 2^pK^p(t)t^{p-1}\int_0^t\E\bigg(\sup\limits_{u\in[0,r]}|X^{n-1}_u-X^{n-2}_u|^p\bigg)\d r.
\end{align}
As for the term $I_2(t)$, owing to $p>1/H$ and $p\geq\theta$, \eqref{MI} and (H) we have
\beg{align}\label{Pf6-The1}
I_2(t)&=\E\left(\sup\limits_{s\in[0,t]}\left|\int_0^s\left(\si(r,\sL_{X^{n-1}_r})-\si(r,\sL_{X^{n-2}_r})\right)\d B_r^H\right|^p\right)\cr
&\leq C_{p,H}t^{pH-1}\int_0^t\|\si(r,\sL_{X^{n-1}_r})-\si(r,\sL_{X^{n-2}_r})\|^p\d r\cr
&\leq C_{p,H}t^{pH-1}\int_0^tK^p(r)\W_\theta(\sL_{X^{n-1}_r},\sL_{X^{n-2}_r})^p\d r\cr
&\leq C_{p,H}K^p(t)t^{pH-1}\int_0^t\E|X^{n-1}_r-X^{n-2}_r|^p\d r\cr
&\leq C_{p,H}K^p(t)t^{pH-1}\int_0^t\E\bigg(\sup\limits_{u\in[0,r]}|X^{n-1}_u-X^{n-2}_u|^p\bigg)\d r.
\end{align}
Plugging \eqref{Pf5-The1} and \eqref{Pf6-The1} into \eqref{Pf4-The1} yields
\beg{align}\label{Pf7-The1}
\E\bigg(\sup\limits_{s\in[0,t]}|X^n_s-X^{n-1}_s|^p\bigg)&\leq2^{p-1}K^p(t)(2^p t^{p-1}+C_{p,H}t^{pH-1})\int_0^t\E\bigg(\sup\limits_{u\in[0,r]}|X^{n-1}_u-X^{n-2}_u|^p\bigg)\d r\nonumber\\
&\leq C_{p,T,H}\int_0^t\E\bigg(\sup\limits_{u\in[0,r]}|X^{n-1}_u-X^{n-2}_u|^p\bigg)\d r
\end{align}
with $C_{p,T,H}:=2^{p-1}K^p(T)(2^p T^{p-1}+C_{p,H}T^{pH-1})$.\\
Hence, by the iteration we arrive at
 \beg{align*}
\E\bigg(\sup\limits_{s\in[0,t]}|X^n_s-X^{n-1}_s|^p\bigg) \leq C_1C_{p,T,H}^n\ff {t^{n-1}} {(n-1)!},
\end{align*}
where $C_1:=\E(\sup_{t\in[0,T]}|X^1_t-\xi|^p)<\infty$ due to Step 1.\\
Consequently, $(X^n)_{n\geq 1}$ is a Cauchy sequence in $\cS^p([0,T])$ with any $p\geq\theta$ and $p>\ff 1 H$,
and then the limit, denoted by $X$, is a solution of \eqref{XX}.

\textit{Step 3. Uniqueness.}
Let $X$ and $Y$ be two solutions of \eqref{XX}.
Along the same lines with Step 2, we derive that as in \eqref{Pf7-The1},
\beg{align*}
\E\bigg(\sup\limits_{s\in[0,t]}|X_s-Y_s|^p\bigg)\leq C_{p,T,H}\int_0^t\E\bigg(\sup\limits_{u\in[0,r]}|X_u-Y_u|^p\bigg)\d r,\ \ t\in[0,T].
\end{align*}
Then, the Gronwall lemma implies that $X_t=Y_t,\ t\in[0,T],\ \P$-a.s.
The proof is now complete.
\qed

\beg{rem}\label{WP-rem}
In \cite[Theorem 3.7]{BJ17}, the authors considered the equation \eqref{XX} with $d=1$, i.e. one-dimensional case.
Under the $\W_1$-Lipschitz conditions and $\xi\in L^2(\Omega\ra\R,\sF_0,\P)$ (namely $\th=1$ and $p=2$ in our result Theorem \ref{(EX)The1} above), they proved the existence and uniqueness of solution $X\in\cS^2([0,T])$ with $\E(\sup_{t\in[0,T]}|X_t|^2)<\infty$ replaced by $\sup_{t\in[0,T]}\E|X_t|^2<\infty$.
So in this sense our result extends and improves \cite[Theorem 3.7]{BJ17}.
\end{rem}




\section{Bismut formulas for the $L$-derivative }

In this section, we consider \eqref{XX} with distribution independent $\si(t)$, that is,
\beg{align}\label{(LDF)Eq1}
\d X_t=b(t,X_t,\sL_{X_t})\d t+\si(t)\d B_t^H,
\end{align}
where $X_0\in L^2(\Omega\ra\R^d,\sF_0,\P)$ with $\sL_{X_0}=\mu$.
The main objective of this section concerns the problem of Bismut formulas for the $L$-derivative of \eqref{(LDF)Eq1}.
We first make necessary preparations concerning the partial derivative with respect to the initial value and the  Malliavin derivative of \eqref{(LDF)Eq1}.
In the second part, we will give a general result about Bismut formula for the $L$-derivative of \eqref{(LDF)Eq1},
the applications to non-degenerate and degenerate cases of \eqref{(LDF)Eq1} are addressed respectively in the last two parts.

\subsection{The partial derivative and the Malliavin derivative of \eqref{(LDF)Eq1}}

We begin with the following assumption.
\beg{enumerate}
\item[(A1)] For every $t\in[0,T]$, $b(t,\cdot,\cdot)\in C^{1,(1,0)}(\R^d\times\sP_2(\R^d))$.
Moreover, there exists a constant $K>0$ such that
\beg{align*}
\|\na b(t,\cdot,\mu)(x)\|+|D^Lb(t,x,\cdot)(\mu)(y)|\le K,\ \ t\in[0,T],\ x,y\in\R^d,\ \mu\in\sP_2(\R^d)
\end{align*}
and $\sup_{t\in[0,T]}(|b(t,0,\de_0)|+\|\si(t)\|)\leq K$.
\end{enumerate}
Note that by the fundamental theorem for Bochner integral (see, for instance, \cite[Proposition A.2.3]{LWbook}) and the definitions of $L$-derivative and the Wasserstein distance, (A1) implies
\beg{align*}
|b(t,x,\mu)-b(t,y,\nu)|\le K(|x-y|+\W_2(\mu,\nu)),\ \  t\in[0,T],\ x,y\in\R^d,\ \mu,\nu\in\sP_2(\R^d).
\end{align*}
Then, it follows from Theorem \ref{(EX)The1} that \eqref{(LDF)Eq1} has a unique solution.

To investigate the partial derivative with respect to initial value of \eqref{(LDF)Eq1}, we first introduce a family of auxiliary equations.
For any $\ve>0$ and $\eta\in L^2(\Om\ra\R^d,\sF_0,\P)$, let $(X^\ve_t)_{t\in[0,T]}$ solve
\beg{align}\label{(LDF)Eq2}
\d X_t^\ve=b(t,X_t^\ve,\sL_{X_t^\ve})\d t+\si(t)\d B_t^H,\ \ X^\ve_0=X_0+\ve\eta,
\end{align}
and define
\beg{align*}
\Upsilon^\ve_t:=\ff{X^\ve_t-X_t}\ve,\ \ t\in[0,T], \ \ \ve>0.
\end{align*}

\beg{lem}\label{UI}
Assume that (A1) holds. Then
\beg{align}\label{UI-1}
\sup_{\ve\in(0,1]}\E\bigg(\sup_{t\in[0,T]}|\Upsilon^\ve_t|^2\bigg)\leq2\e^{8(KT)^2}\E|\eta|^2,
\end{align}
and
\beg{align}\label{Add-UI}
\sup_{\ve\in(0,1],t\in[0,T]}|\Upsilon^\ve_t|^2\leq\left(2|\eta|^2+8(KT)^2\e^{8(KT)^2}\E|\eta|^2\right)\e^{4(KT)^2}.
\end{align}
\end{lem}
\beg{proof}
By \eqref{(LDF)Eq1}-\eqref{(LDF)Eq2} and (A1), we have for any $t\in[0,T]$ and $\ve\in(0,1]$,
\beg{align*}
&\sup\limits_{s\in[0,t]}|X^\ve_s-X_s|^2\cr
&\leq
2\ve^2|\eta|^2+
2\sup\limits_{s\in[0,t]}\left|\int_0^s\left(b(r,X^\ve_r,\sL_{X^\ve_r})-b(r,X_r,\sL_{X_r})\right)\d r\right|^2\cr
&\leq
2\ve^2|\eta|^2+4K^2T\sup\limits_{s\in[0,t]}\int_0^s\left(|X^\ve_r-X_r|^2+\W_2(\sL_{X^\ve_r},\sL_{X_r})^2\right)\d r.
\end{align*}
Taking the expectation on both sides of the above inequality, we get
\beg{align*}
\E\bigg(\sup\limits_{s\in[0,t]}|X^\ve_s-X_s|^2\bigg)\leq2\ve^2\E|\eta|^2+8K^2T\int_0^t\E\bigg(\sup\limits_{u\in[0,r]}|X^\ve_u-X_u|^2\bigg)\d r,
\end{align*}
which implies \eqref{UI-1} and then \eqref{Add-UI} due to the Gronwall inequality.
\end{proof}

With Lemma \ref{UI} in hand, we can present the partial derivative in initial value of the equation  \eqref{(LDF)Eq1}.
Consider now the following linear random ODE on $\R^d$: for any $\eta\in L^2(\Om\ra\R^d,\sF_0,\P)$ and $t\in[0,T]$,
\beg{align}\label{ParV-0}
\d\Gamma^\eta_t&=\left[\na_{\Gamma^\eta_t}b(t,\cdot,\sL_{X_t})(X_t)+\left(\E\<D^Lb(t,y,\cdot)(\sL_{X_t})(X_t), \Gamma^\eta_t\>\right)|_{y=X_t}\right]\d t,\ \ \Gamma^\eta_0=\eta,
\end{align}
where
\beg{align*}
&\E\<D^Lb(t,y,\cdot)(\sL_{X_t})(X_t), \Gamma^\eta_t\>:=\left(\E\<D^Lb_i(t,y,\cdot)(\sL_{X_t})(X_t), \Gamma^\eta_t\>\right)_{1\leq i\leq d}\in\R^d.
\end{align*}
Obviously, (A1) implies that the ODE has a unique solution $\{\Gamma^\eta_t\}_{t\in[0,T]}$ satisfying
\beg{align}\label{ParV-2}
\E\left(\sup_{t\in[0,T]}|\Gamma^\eta_t|^2\right)\leq C_{T,K}\E|\eta|^2.
\end{align}

\beg{prp}\label{Prop-ParV}
Assume that (A1) holds.
Then for any $\eta\in L^2(\Om\ra\R^d,\sF_0,\P)$,
the limit $\na_\eta X_t:=\lim_{\ve\downarrow0}\Upsilon^\ve_t,t\in[0,T]$ exists in $L^2(\Omega\ra C([0,T];\R^d),\P)$
such that $\na_\eta X_t=\Gamma^\eta_t$ holds for each $t\in[0,T]$, i.e., $\na_\eta X_t$ is the unique solution of \eqref{ParV-0}.
\end{prp}

\beg{proof}
To simplify the notation, we denote $X^\ve_\th(t)=X_t+\th(X_t^\ve-X_t), \th\in[0,1]$.
By \eqref{(LDF)Eq1} and \eqref{(LDF)Eq2}, we obtain that for any $t\in[0,T]$,
\beg{align*}
\d\Upsilon^\ve_t&=\ff {b(t,X_t^\ve,\sL_{X_t^\ve})-b(t,X_t,\sL_{X_t})} \ve\d t\cr
&=\left[\ff 1 \ve\int_0^1 \ff \d {\d\th}b(t,X^\ve_\th(t),\sL_{X_t^\ve})\d\th
+\ff 1 \ve\int_0^1 \ff \d {\d\th}b(t,X_t,\sL_{X^\ve_\th(t)})\d\th\right]\d t\cr
&=\left[\int_0^1 \na_{\Upsilon^\ve_t} b(t,\cdot,\sL_{X_t^\ve})(X^\ve_\th(t))\d\th
+\int_0^1 (\E\<D^Lb(t,y,\cdot)(\sL_{X^\ve_\th(t)})(X^\ve_\th(t)),\Upsilon^\ve_t\>)|_{y=X_t}\d\th\right]\d t
\end{align*}
with $\Upsilon^\ve_0=\eta$.
Here, we have used Lemma \ref{FoLD} in the last equality.\\
Then, combining this with \eqref{ParV-0} yields that for each $t\in[0,T]$,
\beg{align*}
\d(\Upsilon^\ve_t-\Gamma^\eta_t)&
=\left[\Phi^\ve_1(t)+\na_{\Upsilon^\ve_t-\Gamma^\eta_t}b(t,\cdot,\sL_{X_t})(X_t)\right]\d t\cr
&\quad+\left[\Phi^\ve_2(t)+\left(\E\<D^Lb(t,y,\cdot)(\sL_{X_t})(X_t),\Upsilon^\ve_t-\Gamma^\eta_t\>\right)|_{y=X_t}\right]\d t, \ \ \Upsilon^\ve_0-\Gamma^\eta_0=0,
\end{align*}
where
\beg{align*}
\Phi^\ve_1(t)&:=\int_0^1\left[\na_{\Upsilon^\ve_t} b(t,\cdot,\sL_{X_t^\ve})(X^\ve_\th(t))-\na_{\Upsilon^\ve_t}b(t,\cdot,\sL_{X_t})(X_t)\right]\d\th,\cr
\Phi^\ve_2(t)&:=\int_0^1\left(\E\<D^Lb(t,y,\cdot)(\sL_{X^\ve_\th(t)})(X^\ve_\th(t))-D^Lb(t,y,\cdot)(\sL_{X_t})(X_t),\Upsilon^\ve_t\>\right)|_{y=X_t}\d\th.
\end{align*}
Consequently, by (A1) we get
\beg{align*}
|\Upsilon^\ve_t-\Gamma^\eta_t|^2\le 4T\int_0^t\left(|\Phi^\ve_1(s)|^2+|\Phi^\ve_2(s)|^2\right)\d s
+4K^2T\int_0^t\left(|\Upsilon^\ve_s-\Gamma^\eta_s|^2+\E|\Upsilon^\ve_s-\Gamma^\eta_s|^2\right)\d s.
\end{align*}
Taking into account of \eqref{UI-1} and \eqref{ParV-2}, the Gronwall inequality leads to
\beg{align}\label{Pf1-Pr-PV}
\E\left(\sup_{t\in[0,T]}|\Upsilon^\ve_t-\Gamma^\eta_t|^2\right)\le4T\e^{8(KT)^2}\int_0^T\E\left(|\Phi^\ve_1(s)|^2+|\Phi^\ve_2(s)|^2\right)\d s.
\end{align}
By the H\"{o}lder inequality and \eqref{UI-1}, one can see that
\beg{align}\label{Pf2-Pr-PV}
&|\Phi^\ve_1(s)|^2+|\Phi^\ve_2(s)|^2\nonumber\\
&\le \int_0^1|\na b(s,\cdot,\sL_{X_s^\ve})(X^\ve_\th(s))-\na b(s,\cdot,\sL_{X_s})(X_s)|^2\d\th\cdot|\Upsilon^\ve_s|^2\nonumber\\
&\quad+\int_0^1\left(\E|D^Lb(s,y,\cdot)(\sL_{X^\ve_\th(s)})(X^\ve_\th(s))-D^Lb(s,y,\cdot)(\sL_{X_s})(X_s)|^2\right)|_{y=X_s}\d\th\cdot
\E|\Upsilon^\ve_s|^2
\end{align}
and
\beg{align*}
\lim_{\ve\downarrow0}\E\bigg(\sup_{\th\in[0,1]}|X^\ve_\th(s)-X_s|^2\bigg)\leq\lim_{\ve\downarrow0}\E|X^\ve_s-X_s|^2=0.
\end{align*}
Then using the condition $b(s,\cdot,\cdot)\in C^{1,(1,0)}(\R^d\times\sP_2(\R^d))$ of (A1) and \eqref{UI-1} again,
we obtain that $|\Phi^\ve_1(s)|^2+|\Phi^\ve_2(s)|^2$ converges to $0$ in probability as $\ve$ goes to $0$.
Additionally, due to \eqref{Pf2-Pr-PV} one has
\beg{align*}
|\Phi^\ve_1(s)|^2+|\Phi^\ve_2(s)|^2\le 4K^2(|\Upsilon^\ve_s|^2+\E|\Upsilon^\ve_s|^2).
\end{align*}
By the dominated convergence theorem and  Lemma \ref{UI}, we conclude that
\beg{align*}
\lim_{\ve\downarrow0}\E\bigg[\sup_{s\in[0,T]}\left(|\Phi^\ve_1(s)|^2+|\Phi^\ve_2(s)|^2\right)\bigg]=0.
\end{align*}
This, along with \eqref{Pf1-Pr-PV}, implies
\beg{align*}
\lim_{\ve\downarrow0}\E\bigg(\sup_{t\in[0,T]}|\Upsilon^\ve_t-\Gamma^\eta_t|^2\bigg)=0,
\end{align*}
which completes the proof.
\end{proof}

For the Malliavin derivative of the equation \eqref{(LDF)Eq1},
consider for each $h\in\cH$ and $\ve>0$ the SDE: for $t\in[0,T]$,
\beg{align}\label{Mall-2}
\d X^{\ve,h}_t=b(t,X_t^{\ve,h},\sL_{X_t})\d t+\si(t)\d(B_t^H+\ve(R_H h)(t)),\ \ X^{\ve,h}_0=X_0.
\end{align}
It is easy to see that under (A1) there exists a unique solution $X^{\ve,h}$ to \eqref{Mall-2}.
Using the pathwise uniqueness of \eqref{(LDF)Eq1} and the fact that $X_t$ can be regarded as a functional of $B^H$ and $X_0$,
the Malliavin directional derivative of $X_t$ along $R_Hh$ is shown by
\beg{align*}
\lim_{\ve\downarrow0}\ff {X^{\ve,h}_t-X_t}\ve
\end{align*}
if the limit exists in $L^2(\Omega\ra C([0,T];\R^d),\P)$.
The above step is partially borrowed from \cite[Proposition 3.5, Page 4762]{RW}.
Noting that $\sL_{X_t}$ in \eqref{Mall-2} is independent of $\ve$,
by the same arguments as in \cite[Lemma 3.1 and Proposition 3.1]{FR} we have the following result.
\beg{prp}\label{Prop-Mall}
Assume that (A1) holds.
Then for any $\eta\in L^2(\Om\ra\R^d,\sF_0,\P)$ and $h\in\cH$, the limit
\beg{align*}
\D_{R_Hh}X_t:=\lim_{\ve\downarrow0}\ff {X^{\ve,h}_t-X_t}\ve, \ \ t\in[0,T]
\end{align*}
exists in $L^2(\Omega\ra C([0,T];\R^d),\P)$ such that $\D_{R_Hh}X_t=(\<\D X^i_t,h\>_\cH)_{1\leq i\leq d}\in\R^d$ holds for every $t\in[0,T]$ and satisfies
\beg{align}\label{Mall-1}
\D_{R_Hh}X_t=\int_0^t\na_{\D_{R_Hh}X_s}b(s,\cdot,\sL_{X_s})(X_s)\d s+\int_0^t\si(s)\d(R_H h)(s), \ \ t\in[0,T].
\end{align}
\end{prp}

\subsection{Bismut formula: a general result}

In this part, we aim to establish a general result of Bismut formula for the $L$-derivative of \eqref{(LDF)Eq1},
which is then applied to non-degenerate and degenerate cases in the next two parts, respectively.
More precisely, for any $\mu\in\sP_2(\R^d)$, let $(X_t^\mu)_{t\in[0,T]}$ be the solution to \eqref{(LDF)Eq1} with $\sL_{X_0}=\mu$
and denote $P^*_t\mu=\sL_{X_t^\mu}$ for every $t\in[0,T]$.
Now, define
\beg{align*}
(P_t f)(\mu):=\int_{\R^d} f\d(P^*_t\mu)=\E f(X_t^\mu), \ \ t\in[0,T], f\in\sB_b(\R^d), \mu\in\sP_2(\R^d).
\end{align*}
For any $t\in(0,T],\mu\in\sP_2(\R^d)$ and $\phi\in L^2(\R^d\ra\R^d,\mu)$, we are to find an integrable random variable $M_t(\mu,\phi)$ such that
\beg{align*}
D^L_\phi(P_t f)(\mu)=\E\left(f(X_t^\mu)M_t(\mu,\phi)\right),\ \ f\in\sB_b(\R^d).
\end{align*}
Toward this goal, for any $\ve\in[0,1]$ and  $\phi\in L^2(\R^d\ra\R^d,\mu)$, let $X_t^{\mu_{\varepsilon,\phi}}$ denote the solution of \eqref{(LDF)Eq1} with $X_0^{\mu_{\varepsilon,\phi}}=(\mathrm{Id}+\varepsilon\phi)(X_0)$.
According to Proposition \ref{Prop-ParV} and \ref{Prop-Mall},
$\na_{\phi(X_0)} X_\cdot^{\mu_{\ve,\phi}}$ and $\D_{R_H h_{s_0}^{\ve,\phi}}X_\cdot^{\mu_{\ve,\phi}}$ below are both well-defined with any $s_0\in[0,T)$,
and moreover satisfy \eqref{ParV-0} with $\eta=\phi(X_0)$ and \eqref{Mall-1}, respectively.
In order to ease notations, we simply write $\mu_{\varepsilon,\phi}=\sL_{(\mathrm{Id}+\varepsilon\phi)(X_0)}$,
and if $s_0=0$ or $\ve=0$, we often suppress $s_0$ or $\ve$ (e.g., $R_H h_0^{\ve,\phi}=R_H h^{\ve,\phi}, h_0^{0,\phi}=h^\phi, X_t^{\mu_{\varepsilon,\phi}}=X_t^\mu,\cdots,$ etc.).

Our main result is the following.

\beg{thm}\label{(Ge)Th1}
Assume that (A1) holds, and that for any  $\ve\in[0,1]$ and $s_0\in[0,T)$, there exists $h_{s_0}^{\ve,\phi}\in\mathrm{Dom}\delta\cap\cH$ such that
$\D_{R_H h_{s_0}^{\ve,\phi}}X_T^{\mu_{\varepsilon,\phi}}=\na_{\phi(X_0)} X_T^{\mu_{\varepsilon,\phi}}$ and $(R_H h_{s_0}^{\ve,\phi})(t)=0$ for all $t\in[0,s_0]$.
Assume in addition that
$\int_0^1\left(\E\delta^2(h_{s_0}^{\tau,\phi})\right)^\ff 1 2\d\tau<\infty$ and
$$\lim_{\ve\ra 0^+}\E|\delta(h^{\ve,\phi})-\delta(h^\phi)|=0,\ \  \phi\in L^2(\R^d\ra\R^d,\mu).$$
Then we have\\
(i) For any $f\in\sB_b(\R^d)$, $P_Tf$ is weakly $L$-differentiable at $\mu$, and moreover
\beg{align*}
D^L_\phi(P_T f)(\mu)=\E(f(X_T^\mu)\delta(h^\phi)),\ \  \phi\in L^2(\R^d\ra\R^d,\mu).
\end{align*}
(ii) For any $f\in\sB_b(\R^d), P_T f$ is $L$-differentiable at $\mu$, if $\E\delta^2(h^\phi)\leq \tilde{L}\|\phi\|^2_{L^2_\mu}$ with a constant $\tilde{L}>0$ and
\beg{align}\label{(Ge)Th1-1}
 \lim_{\|\phi\|_{L^2_\mu}\ra 0}\sup\limits_{\ve\in(0,1]}\ff{\E|\delta(h^{\ve,\phi})-\delta(h^\phi)|}{\|\phi\|_{L^2_\mu}}=0.
\end{align}
\end{thm}

In order to prove the theorem, we first give the following lemma which will play a crucial role in the proof.

\beg{lem}\label{(Ge)lem1}
Assume that (A1) holds, and that for any  $\ve\in[0,1]$ and $s_0\in[0,T)$, there exists $h_{s_0}^{\ve,\phi}\in\mathrm{Dom}\delta\cap\cH$ such that
$\D_{R_H h_{s_0}^{\ve,\phi}}X_T^{\mu_{\varepsilon,\phi}}=\na_{\phi(X_0)} X_T^{\mu_{\varepsilon,\phi}}$ and $(R_H h_{s_0}^{\ve,\phi})(t)=0$ for all $t\in[0,s_0]$.
Then for any $\ve\in[0,1], s_0\in[0,T)$ and $f\in\sB_b(\R^d)$,
\beg{align*}
\E(f(X_T^{\mu_{\ve,\phi}})-f(X_T^\mu)|\sF_{s_0})=\int_0^\ve\E\left(f(X_T^{\mu_{\tau,\phi}})\delta(h_{s_0}^{\tau,\phi})\Big|\sF_{s_0}\right)\d\tau.
\end{align*}
In particular, it holds
\beg{align*}
\E(f(X_T^{\mu_{\ve,\phi}})-f(X_T^\mu))=\int_0^\ve\E\left(f(X_T^{\mu_{\tau,\phi}})\delta(h^{\tau,\phi})\right)\d\tau.
\end{align*}

\end{lem}

\beg{proof}
Since $\D_{R_H h_{s_0}^{\ve,\phi}}X_T^{\mu_{\varepsilon,\phi}}=\na_{\phi(X_0)} X_T^{\mu_{\varepsilon,\phi}}$,
we deduce that for any $f\in C^1_b(\R^d)$,
\beg{align}\label{(Ge)lem1-1}
&\E(f(X_T^{\mu_{\ve,\phi}})-f(X_T^\mu)|\sF_{s_0})
=\E\left(\int_0^\ve\ff \d {\d\tau}f(X_T^{\mu_{\tau,\phi}})\d \tau\Big|\sF_{s_0}\right)\cr
&=\E\left(\int_0^\ve\<\na f(X_T^{\mu_{\tau,\phi}}),\na_{\phi(X_0)} X_T^{\mu_{\tau,\phi}}\>\d\tau\Big|\sF_{s_0}\right)\cr
&=\int_0^\ve\E\left(\<\na f(X_T^{\mu_{\tau,\phi}}),\na_{\phi(X_0)} X_T^{\mu_{\tau,\phi}}\>\Big|\sF_{s_0}\right)\d\tau\cr
&=\int_0^\ve\E\left(\<\na f(X_T^{\mu_{\tau,\phi}}),
\D_{R_H h_{s_0}^{\tau,\phi}}X_T^{\mu_{\tau,\phi}}\>\Big|\sF_{s_0}\right)\d\tau\cr
&=\int_0^\ve\E\left(\D_{R_H h_{s_0}^{\tau,\phi}}f(X_T^{\mu_{\tau,\phi}})
\Big|\sF_{s_0}\right)\d\tau\cr
&=\int_0^\ve\E\left(\<\D f(X_T^{\mu_{\tau,\phi}}),h_{s_0}^{\tau,\phi}\>_\cH
\Big|\sF_{s_0}\right)\d\tau.
\end{align}
Now, let $\zeta\in\mathrm{Dom}\D$ be any bounded and $\sF_{s_0}$-measurable smooth random variable,
by \cite[Proposition 1.2.3]{ND} we have for any $\tau\in[0,\ve]$,
\beg{align}\label{(Ge)lem1-2}
&\E\left(\zeta\<\D f(X_T^{\mu_{\tau,\phi}}),h_{s_0}^{\tau,\phi}\>_\cH\right)\cr
&=\E\left[\left\<\D(\zeta f(X_T^{\mu_{\tau,\phi}})),h_{s_0}^{\tau,\phi}\right\>_\cH-f(X_T^{\mu_{\tau,\phi}})\<\D \zeta,h_{s_0}^{\tau,\phi}\>_\cH\right]\cr
&=\E\left[\zeta f(X_T^{\mu_{\tau,\phi}})\delta(h_{s_0}^{\tau,\phi})-f(X_T^{\mu_{\tau,\phi}})\<\D \zeta,h_{s_0}^{\tau,\phi}\>_\cH\right]\cr
&=\E\left[\zeta f(X_T^{\mu_{\tau,\phi}})\delta(h_{s_0}^{\tau,\phi})-f(X_T^{\mu_{\tau,\phi}})\<K_H^*\D \zeta,K_H^*h_{s_0}^{\tau,\phi}\>_{L^2([0,T];\R^d)}\right],
\end{align}
where the last equality is due to the fact that $K_H^*$ is an isometry between $\cH$ and a closed subspace of $L^2([0,T];\R^d)$.\\
Using Proposition \ref{Tra-Prp1} and the fact that $(\D^W\zeta)(t)=0$ for all $t>s_0$, we get
\beg{align}\label{(Ge)lem1-3}
&\left\<K_H^*\D \zeta,K_H^*h_{s_0}^{\tau,\phi}\right\>_{L^2([0,T];\R^d)}=\left\<\D^W \zeta,K_H^*h_{s_0}^{\tau,\phi}\right\>_{L^2([0,T];\R^d)}\cr
&=\int_0^T\left\<(\D^W\zeta)(t),(K_H^*h_{s_0}^{\tau,\phi})(t)\right\>\d t=\int_0^{s_0}\left\<(\D^W\zeta)(t),(K_H^*h_{s_0}^{\tau,\phi})(t)\right\>\d t=0.
\end{align}
Here we have used  $K_H^*h_{s_0}^{\tau,\phi}=K_H^{-1}(R_Hh_{s_0}^{\tau,\phi})$ and the fact that $(R_H h_{s_0}^{\tau,\phi})(t)=0$ for $t\in[0,s_0]$
in the last equality.\\
Substituting \eqref{(Ge)lem1-3} into \eqref{(Ge)lem1-2} implies
\beg{align*}
\E\left[\zeta\<\D f(X_T^{\mu_{\tau,\phi}}),h_{s_0}^{\tau,\phi}\>_\cH\right]
&=\E\left[\zeta f(X_T^{\mu_{\tau,\phi}})\delta(h_{s_0}^{\tau,\phi})\right].
\end{align*}
Hence, combining this with \eqref{(Ge)lem1-1} we obtain
\beg{align}\label{(Ge)lem1-4}
\E(f(X_T^{\mu_{\ve,\phi}})-f(X_T^\mu)|\sF_{s_0})
=\int_0^\ve\E\left(f(X_T^{\mu_{\tau,\phi}})\delta(h_{s_0}^{\tau,\phi})\Big|\sF_{s_0}\right)\d\tau,\ \ f\in C^1_b(\R^d).
\end{align}
Set
\beg{align*}
\nu^{\ve,\phi}_{s_0}(A):=\int_0^\ve\E\left(
\mathrm{I}_A (X_T^{\mu_{\tau,\phi}})|\delta(h_{s_0}^{\tau,\phi})|\right)\d\tau,\ \ A\in\sB(\R^d),
\end{align*}
which is a finite measure on $\R^d$.
Then $C^1_b(\R^d)$ is dense in $L^1(\R^d,\sL_{X_T^{\mu_{\ve,\phi}}}+\sL_{X_T^\mu}+\nu^{\ve,\phi}_{s_0})$.
Therefore, \eqref{(Ge)lem1-4} holds for any $f\in\sB_b(\R^d)$.
This completes the proof.
\end{proof}

\emph{Proof of Theorem \ref{(Ge)Th1}.}
We divide the proof into two steps.

\textsl{Step 1. Claim: For any $f\in\sB_b(\R^d), P_T f$ is weakly $L$-differentiable at $\mu=\sL_{X_0}$
(namely $(P_T f)(\mu\circ(\mathrm{Id}+\cdot)^{-1}): L^2(\R^d\ra\R^d,\mu)\ra\R$ is G\^{a}teaux differentiable at $0$),
and moreover $D^L_\phi(P_T f)(\mu)=\E(f(X_T^\mu)\delta(h^\phi))$ holds for each $\phi\in L^2(\R^d\ra\R^d,\mu)$.}
Due to Lemma \ref{(Ge)lem1}, we deduce that for any $f\in\sB_b(\R^d)$ and $\phi\in L^2(\R^d\ra\R^d,\mu)$,
\beg{align}\label{Pf1(Ge)Th1}
&\ff{(P_T f)(\mu\circ(\mathrm{Id}+\ve\phi)^{-1})-(P_T f)(\mu)}{\ve}-\E(f(X_T^\mu)\delta(h^\phi))\cr
&=\ff{(P_T f)(\sL_{(\mathrm{Id}+\varepsilon\phi)(X_0)})-(P_T f)(\sL_{X_0})}{\varepsilon}-\E(f(X_T^\mu)\delta(h^\phi))\cr
&=\ff{\E f(X_T^{\mu_{\ve,\phi}})-\E f(X_T^\mu)}{\ve}-\E(f(X_T^\mu)\delta(h^\phi))\cr
&=\frac{1}{\ve}\int_0^\ve\E\left(f(X_T^{\mu_{\tau,\phi}})\delta(h^{\tau,\phi})\right)\d\tau-\E(f(X_T^\mu)\delta(h^\phi))\cr
&=\frac{1}{\ve}\int_0^\ve\E\left[f(X_T^{\mu_{\tau,\phi}})(\delta(h^{\tau,\phi})-\delta(h^\phi))\right]\d\tau
+\frac{1}{\ve}\int_0^\ve\E\left[(f(X_T^{\mu_{\tau,\phi}})-f(X_T^\mu))\delta(h^\phi)\right]\d\tau\cr
&=:I_1(\phi)+I_2(\phi).
\end{align}
Since $\lim_{\tau\ra 0^+}\E|\delta(h^{\tau,\phi})-\delta(h^\phi)|=0$ for any $\phi\in L^2(\R^d\ra\R^d,\mu)$, we obtain
\beg{align}\label{Pf4(Ge)Th1}
\limsup_{\ve\ra 0^+}|I_1(\phi)|&\leq\|f\|_\infty\lim_{\ve\ra 0^+} \frac{1}{\ve}\int_0^\ve\E|\delta(h^{\tau,\phi})-\delta(h^\phi)|\d\tau\cr
&=\|f\|_\infty\lim_{\ve\ra 0^+}\E|\delta(h^{\ve,\phi})-\delta(h^\phi)|=0.
\end{align}
For $I_2(\phi)$, we get for any $s_0\in(0,T)$,
\beg{align}\label{Pf2(Ge)Th1}
|I_2(\phi)|&\leq\frac{1}{\ve}\int_0^\ve\left|\E\left[(f(X_T^{\mu_{\tau,\phi}})-f(X_T^\mu))(\delta(h^\phi)-\E(\delta(h^\phi)|\sF_{s_0}))\right]\right|\d\tau\cr
&\quad+\frac{1}{\ve}\int_0^\ve\left|\E\left[(f(X_T^{\mu_{\tau,\phi}})-f(X_T^\mu))\E(\delta(h^\phi)|\sF_{s_0})\right]\right|\d\tau\cr
&\leq2\|f\|_\infty  \E\left|\delta(h^\phi)-\E(\delta(h^\phi)|\sF_{s_0})\right|\cr
&\quad+\frac{1}{\ve}\int_0^\ve\left|\E\left[(f(X_T^{\mu_{\tau,\phi}})-f(X_T^\mu))\E(\delta(h^\phi)|\sF_{s_0})\right]\right|\d\tau. \end{align}
On one hand, it is easy to see that
\beg{align}\label{Pf3(Ge)Th1}
\lim_{s_0\ra T^-}\lim_{\ve\ra 0^+}\E\left|\delta(h^\phi)-\E(\delta(h^\phi)|\sF_{s_0})\right|
=\lim_{s_0\ra T^-}\E\left|\delta(h^\phi)-\E(\delta(h^\phi)|\sF_{s_0})\right|=0.
\end{align}
On the other hand, note that by Lemma \ref{(Ge)lem1} again, we have
\beg{align*}
&|\E[(f(X_T^{\mu_{\tau,\phi}})-f(X_T^\mu))\E(\delta(h^\phi)|\sF_{s_0})]|\cr
&=|\E[\E(\delta(h^\phi)|\sF_{s_0})\E(f(X_T^{\mu_{\tau,\phi}})-f(X_T^\mu)|\sF_{s_0})]|\cr
&=\left|\E\left[\E(\delta(h^\phi)|\sF_{s_0})\int_0^\tau\E\left(f(X_T^{\mu_{\theta,\phi}})\delta(h_{s_0}^{\theta,\phi})\Big|\sF_{s_0}\right)\d\theta\right]\right|\cr
&=\left|\int_0^\tau\E\left[\E(\delta(h^\phi)|\sF_{s_0})f(X_T^{\mu_{\theta,\phi}})\delta(h_{s_0}^{\theta,\phi})\right]\d\theta\right|\cr
&\leq\|f\|_\infty(\E\delta^2(h^\phi))^\ff 1 2\int_0^\tau\left(\E\delta^2(h_{s_0}^{\theta,\phi})\right)^\ff 1 2\d\theta,
\end{align*}
which goes to zero as $\tau\ra0$ because of $\int_0^1\left(\E\delta^2(h_{s_0}^{\theta,\phi})\right)^\ff 1 2\d\theta<\infty$.
This means that the function $\tau\mapsto\E[(f(X_T^{\mu_{\tau,\phi}})-f(X_T^\mu))\E(\delta(h^\phi)|\sF_{s_0})]$ is continuous at 0.
Then, we derive that for each $s_0\in(0,T)$,
\beg{align}\label{Pf5(Ge)Th1}
\lim_{\ve\ra 0^+}\frac{1}{\ve}\int_0^\ve\left|\E\left[(f(X_T^{\mu_{\tau,\phi}})-f(X_T^\mu))\E(\delta(h^\phi)|\sF_{s_0})\right]\right|\d\tau
=0.
\end{align}
Hence, plugging \eqref{Pf3(Ge)Th1} and \eqref{Pf5(Ge)Th1} into \eqref{Pf2(Ge)Th1} implies that $\lim_{\ve\ra 0^+}|I_2(\phi)|=0$.
Combining this and \eqref{Pf4(Ge)Th1} with \eqref{Pf1(Ge)Th1} yields the desired assertion.

\textsl{Step 2. Claim: For any $f\in\sB_b(\R^d), P_T f$ is $L$-differentiable at $\mu=\sL_{X_0}$ (namely $(P_T f)(\mu\circ(\mathrm{Id}+\cdot)^{-1}): L^2(\R^d\ra\R^d,\mu)\ra\R$ is Fr\'{e}chet differentiable at $0$).}
According to the definition of the $L$-derivative, it is sufficient to show that for any $f\in\sB_b(\R^d)$,
\beg{align*}
\lim_{\|\phi\|_{L^2_\mu}\ra 0}\ff{(P_T f)(\mu\circ(\mathrm{Id}+\phi)^{-1})-(P_T f)(\mu)-\E(f(X_T^\mu)\delta(h^\phi))}{\|\phi\|_{L^2_\mu}}=0.
\end{align*}
Applying Lemma \ref{(Ge)lem1} with $\varepsilon=1$, we deduce that for any $f\in\sB_b(\R^d)$,
\beg{align*}
&\ff{|(P_T f)(\mu\circ(\mathrm{Id}+\phi)^{-1})-(P_T f)(\mu)-\E(f(X_T^\mu)\delta(h^\phi))|}{\|\phi\|_{L^2_\mu}}\cr
&=\ff{|\E f(X_T^{\mu_{1,\phi}})-\E f(X_T^\mu) -\E(f(X_T^\mu)\delta(h^\phi))|}{\|\phi\|_{L^2_\mu}}\cr
&=\ff{\left|\int_0^1[\E(f(X_T^{\mu_{\tau,\phi}})\delta(h^{\tau,\phi}))-\E(f(X_T^\mu)\delta(h^\phi))]\d\tau\right|}{\|\phi\|_{L^2_\mu}}\cr
&\leq \ff{\left|\int_0^1\E[f(X_T^{\mu_{\tau,\phi}})(\delta(h^{\tau,\phi})-\delta(h^\phi))]\d\tau\right|}{\|\phi\|_{L^2_\mu}}
+\ff{\left|\int_0^1\E[(f(X_T^{\mu_{\tau,\phi}})-f(X_T^\mu))\delta(h^\phi)]\d\tau\right|}{\|\phi\|_{L^2_\mu}}\cr
&\leq \ff{\|f\|_\infty\int_0^1\E|\delta(h^{\tau,\phi})-\delta(h^\phi)|\d\tau}{\|\phi\|_{L^2_\mu}}
+\tilde{L}\int_0^1(\E|f(X_T^{\mu_{\tau,\phi}})-f(X_T^\mu)|^2)^\ff 1 2\d\tau\cr
&=:J_1(\phi)+J_2(\phi),
\end{align*}
where the last inequality is due to the condition $\E\delta^2(h^\phi)\leq \tilde{L}\|\phi\|^2_{L^2_\mu}$. \\
Obviously, it follows from \eqref{(Ge)Th1-1}  that $\lim_{\|\phi\|_{L^2_\mu}\ra 0}J_1(\phi)=0$.  \\
For $J_2(\phi)$, note first that by the Lusin theorem (see, e.g., \cite[Theorem 7.4.4]{Coh}),
there exist $\{f_n\}_{n\geq1}\subset C_b(\R^d)$ and compact sets $\{K_n\}_{n\geq1}$ such that
\beg{align*}
 f_n|_{K_n}=f|_{K_n}, \ \ \ \|f_n\|_\infty\leq\|f\|_\infty, \ \ \ (\sL_{X_T^{\mu_{\tau,\phi}}}+\sL_{X_T^\mu})(K_n^c)\leq\ff 1 {n^2}.
\end{align*}
Then, we obtain
\beg{align}\label{Pf6(Ge)Th1}
&(\E|f(X_T^{\mu_{\tau,\phi}})-f(X_T^\mu)|^2)^\ff 1 2\cr
&\leq(\E|f(X_T^{\mu_{\tau,\phi}})-f_n(X_T^{\mu_{\tau,\phi}})|^2)^\ff 1 2+(\E|f_n(X_T^{\mu_{\tau,\phi}})-f_n(X_T^\mu)|^2)^\ff 1 2\cr
&\quad+(\E|f_n(X_T^\mu)-f(X_T^\mu)|^2)^\ff 1 2\cr
&\leq \ff {4\|f\|_\infty} n+(\E|f_n(X_T^{\mu_{\tau,\phi}})-f_n(X_T^\mu)|^2)^\ff 1 2.   \end{align}
Note that for any $\tau\in[0,1]$, we have
\beg{align*}
\limsup_{\|\phi\|_{L^2_\mu}\ra 0}\E|X_T^{\mu_{\tau,\phi}}-X_T^\mu|^2\leq C_{T,K}\lim_{\|\phi\|_{L^2_\mu}\ra 0}\|\phi\|^2_{L^2_\mu}=0.
\end{align*}
Consequently, the dominated convergence theorem yields that for every $n\geq 1$,
\beg{align*}
\lim_{\|\phi\|_{L^2_\mu}\ra 0}\E|f_n(X_T^{\mu_{\tau,\phi}})-f_n(X_T^\mu)|^2=0.
\end{align*}
Combining this with \eqref{Pf6(Ge)Th1} yields
\beg{align*}
\lim_{\|\phi\|_{L^2_\mu}\ra 0}(\E|f(X_T^{\mu_{\tau,\phi}})-f(X_T^\mu)|^2)^\ff 1 2=0.
\end{align*}
By the dominated convergence theorem again, we obtain that $\lim_{\|\phi\|_{L^2_\mu}\ra 0}J_2(\phi)=0$,
which completes the proof.
\qed

\subsection{Bismut formula: the non-degenerate case}

This part is devoted to applying our general Theorem \ref{(Ge)Th1} to the non-degenerate case of \eqref{(LDF)Eq1}.
In additional to (A1), we also need the following assumptions.
\beg{enumerate}
\item[(A2)] There exists a constant $\tilde{K}>0$ such that

\item[(i)]  for any $t,s\in[0,T],\ x,y,z_1,z_2\in\R^d,\ \mu,\nu\in\sP_2(\R^d)$,
\beg{align*}
&\|\na b(t,\cdot,\mu)(x)-\na b(s,\cdot,\nu)(y)\|+|D^Lb(t,x,\cdot)(\mu)(z_1)-D^Lb(s,y,\cdot)(\nu)(z_2)|\cr
&\le \tilde{K}(|t-s|^{\alpha_0}+|x-y|^{\beta_0}+|z_1-z_2|^{\gamma_0}+\W_2(\mu,\nu)),
\end{align*}
where $\alpha_0\in(H-1/2,1]$ and $\beta_0,\gamma_0\in(1-1/(2H),1]$.

\item[(ii)] $\si$ is invertible and $\si^{-1}$ is H\"{o}lder continuous of order $\delta_0\in(H-1/2,1]$:
\beg{align*}
\|\si^{-1}(t)-\si^{-1}(s)\|\le\tilde{K}|t-s|^{\delta_0}, \ \  t,s\in[0,T].
\end{align*}

\item[(A3)] The derivatives
$$\partial_t(D^Lb(\cdot,x,\cdot)(\mu)(y))(t),\ \na(D^Lb(t,\cdot,\cdot)(\mu)(y))(x),$$
$$D^L(D^Lb(t,x,\cdot)(\cdot)(y))(\mu)(z),\ \na(D^Lb(t,x,\cdot)(\mu)(\cdot))(y)$$
exist and are bounded continuous in the corresponding arguments $(t,x,\mu,y)$ or $(t,x,\mu,y,z)$.
We denote the bounded constants by a common one $\bar K>0$.
\end{enumerate}

Our main goal in the current part is to prove the following result.

\beg{thm}\label{(ND)Th}
Assume that (A1), (A2) and (A3) hold.
Then for any $\mu\in\sP_2(\R^d)$ and $f\in\sB_b(\R^d), P_T f$ is $L$-differentiable at $\mu$ and
\beg{align*}
D^L_\phi(P_T f)(\mu)=\E\left(f(X_T^\mu)\int_0^T\<K_H^{-1}(R_H h^\phi)(t),\d W_t\>\right),\ \  \phi\in L^2(\R^d\ra\R^d,\mu),
\end{align*}
where $h^\phi\in\mathrm{Dom}\delta\cap\cH$ and satisfies for every $t\in[0,T]$,
\beg{align*}
(R_H h^\phi)(t)
&=\int_0^t\si^{-1}(s)\left[\ff 1 T \na_{\phi(X_0)} X_s^\mu+\ff s T\left(\E\<D^Lb(s,y,\cdot)(\sL_{X_s^\mu})(X_s^\mu), \na_{\phi(X_0)} X_s^\mu\>\right)|_{y=X_s^\mu}\right]\d s.
\end{align*}
\end{thm}

\beg{rem}\label{Non-Re1}
By \eqref{InOp}, one can recast the term $K_H^{-1}(R_H h^\phi)(t)$ in the theorem as
\beg{align*}
&K_H^{-1}(R_H h^\phi)(t)=\frac{(H-\ff 1 2)t^{H-\ff 1 2}}{\Gamma(\frac{3}{2}-H)}
\Bigg[\ff {t^{1-2H}\si^{-1}(t)\varrho(t)}{H-\ff 1 2}
+\si^{-1}(t)\varrho(t)\int_0^t\frac{t^{\frac{1}{2}-H}-s^{\frac{1}{2}-H}}{(t-s)^{\frac{1}{2}+H}}\d s\cr
&\qquad\qquad\qquad\qquad\qquad+\varrho(t)\int_0^t\frac{\si^{-1}(t)-\si^{-1}(s)}{(t-s)^{\frac{1}{2}+H}}s^{\ff 1 2-H}\d s
+\int_0^t\frac{\varrho(t)-\varrho(s)}{(t-s)^{\frac{1}{2}+H}}\si^{-1}(s)s^{\ff 1 2-H}\d s\Bigg],
\end{align*}
where for any $s\in[0,T],\varrho(s)=\ff 1 T \na_{\phi(X_0)} X_s^\mu+\ff s T\left(\E\<D^Lb(s,y,\cdot)(\sL_{X_s^\mu})(X_s^\mu), \na_{\phi(X_0)} X_s^\mu\>\right)|_{y=X_s^\mu}$.
\end{rem}

In order to prove the theorem, we present the following lemma, in which the $L^2$-norm conditionally to $\sF_0$ and the $L^1$-norm
error estimates between $\na_{\phi(X_0)}X_t^{\mu_{\ve,\phi}}$ and $\na_{\phi(X_0)}X_t^\mu$ are provided, respectively.

\beg{lem}\label{(ND)Le}
Assume that (A1) and (A2) are satisfied.
Then for any $t\in[0,T]$,
\beg{align}\label{(ND)Le1}
\E|\na_{\phi(X_0)}X_t^{\mu_{\ve,\phi}}-\na_{\phi(X_0)}X_t^\mu|\leq C_{T,K,\tilde{K}}\ell(\ve,\phi)\|\phi\|_{L^2_\mu}
\end{align}
and
\beg{align}\label{(ND)Le2}
&\E(|\na_{\phi(X_0)}X_t^{\mu_{\ve,\phi}}-\na_{\phi(X_0)}X_t^\mu|^2|\sF_0)\cr
&\le C_{T,K,\tilde{K}}\left(\tilde{\ell}_1^2(\ve,\phi)\|\phi\|_{L^2_\mu}^2+\tilde{\ell}_2^2(\ve,\phi)\|\phi\|_{L^2_\mu}^2+\tilde{\ell}_3^2(\ve,\phi)|\phi(X_0)|^2\right),
\end{align}
where
\beg{align}\label{AD0Le2}
&\ell(\ve,\phi)=\ve^{\beta_0}\|\phi\|_{L^2_\mu}^{\beta_0}+\ve^{\gamma_0}\|\phi\|_{L^2_\mu}^{\gamma_0}+\ve\|\phi\|_{L^2_\mu},\\
&\tilde{\ell}_1(\ve,\phi)
=\ve^{\beta_0}\|\phi\|_{L^2_\mu}^{\beta_0}+\ve^{\gamma_0}\|\phi\|_{L^2_\mu}^{\gamma_0}+\ve\|\phi\|_{L^2_\mu}
+\ve^{\ff {\beta_0} 2}\|\phi\|_{L^2_\mu}^{\ff {\beta_0} 2}+\ve^{\ff 1 2}\|\phi\|^{\ff 1 2}_{L^2_\mu},\label{AD1Le2} \\
&\tilde{\ell}_2(\ve,\phi)
=\ve^{\ff {\beta_0}2}|\phi(X_0)|^{\ff {\beta_0} 2}+\ve^{\beta_0}|\phi(X_0)|^{\beta_0},\label{AD2Le2} \\
&\tilde{\ell}_3(\ve,\phi)
=\ve^{\ff {\beta_0} 2}\|\phi\|_{L^2_\mu}^{\ff {\beta_0} 2}+\ve^{\ff 1 2}\|\phi\|^{\ff 1 2}_{L^2_\mu}+\ve^{\ff {\beta_0}2}|\phi(X_0)|^{\ff {\beta_0} 2}.\label{AD3Le2}
\end{align}
\end{lem}

\beg{rem}\label{(ND)Re1}
By a straightforward calculation, one can see that
\beg{align*}
&\lim_{\ve\ra0}\left[\ell(\ve,\phi)+\tilde{\ell}_1(\ve,\phi)+\E\left(\tilde{\ell}^2_2(\ve,\phi)+\tilde{\ell}^2_3(\ve,\phi)\right)\right]=0,\cr
&\lim_{|\phi\|_{L^2_\mu}\ra0}\sup_{\ve\in(0,1]}\left[\ell(\ve,\phi)+\tilde{\ell}_1(\ve,\phi))
+\E\left(\tilde{\ell}^2_2(\ve,\phi)+\tilde{\ell}^2_3(\ve,\phi)\right)\right]=0.
\end{align*}

\end{rem}

\beg{proof}
By Proposition \ref{Prop-ParV} with $\eta=\phi(X_0)$, we get for every $t\in[0,T]$,
\beg{align*}
\na_{\phi(X_0)}X_t^{\mu_{\ve,\phi}}-\na_{\phi(X_0)}X_t^\mu
&=\int_0^t\bigg[\na_{\na_{\phi(X_0)}X_s^{\mu_{\ve,\phi}}}b(s,\cdot,\sL_{X_s^{\mu_{\ve,\phi}}})(X_s^{\mu_{\ve,\phi}})\cr
&\qquad\quad-\na_{\na_{\phi(X_0)}X_s^\mu} b(s,\cdot,\sL_{X_s^\mu})(X_s^\mu)\cr
&\qquad\quad+\left(\E\<D^Lb(s,y,\cdot)(\sL_{X_s^{\mu_{\ve,\phi}}})(X_s^{\mu_{\ve,\phi}}), \na_{\phi(X_0)} X_s^{\mu_{\ve,\phi}}\>\right)\Big|_{y=X_s^{\mu_{\ve,\phi}}}\cr
&\qquad\quad-\left(\E\<D^Lb(s,y,\cdot)(\sL_{X_s^\mu})(X_s^\mu), \na_{\phi(X_0)} X_s^\mu\>\right)|_{y=X_s^\mu}\bigg]\d s.
\end{align*}
Let $\zeta_t=|\na_{\phi(X_0)}X_t^{\mu_{\ve,\phi}}-\na_{\phi(X_0)}X_t^\mu|$.
Then, by (A1) and (A2) we have for any $t\in[0,T]$,
\beg{align*}
&\E(\zeta_t|\sF_0)\le K\int_0^t(\E(\zeta_s|\sF_0)+\E\zeta_s)\d s\cr
&+\tilde{K}\int_0^t\E\left((|X_s^{\mu_{\ve,\phi}}-X_s^\mu|^{\beta_0}+\W_2(\sL_{X_s^{\mu_{\ve,\phi}}},\sL_{X_s^\mu}))|\na_{\phi(X_0)}X_s^\mu|\big|\sF_0\right)\d s\cr
&+\tilde{K}\int_0^t\left(\E(|X_s^{\mu_{\ve,\phi}}-X_s^\mu|^{\beta_0}|\sF_0)+\W_2(\sL_{X_s^{\mu_{\ve,\phi}}},\sL_{X_s^\mu})\right)\E|\na_{\phi(X_0)}X_s^\mu|\d s\cr
&+\tilde{K}\int_0^t\E\left(|X_s^{\mu_{\ve,\phi}}-X_s^\mu|^{\gamma_0}|\na_{\phi(X_0)}X_s^\mu|\right)\d s.
\end{align*}
Notice that by (A1), we derive for any $p>0$,
\beg{align}\label{Pf2(ND)Le}
\sup_{s\in[0,T]}\E(|X_s^{\mu_{\ve,\phi}}- X_s^\mu|^p|\sF_0)\leq C_{p,T,K}\ve^p\left(\|\phi\|_{L^2_\mu}^p+|\phi(X_0)|^p\right)
\end{align}
and
\beg{align}\label{Pf3(ND)Le}
\sup_{s\in[0,T],\ve\in[0,1]}\E(|\na_{\phi(X_0)} X_s^{\mu_{\ve,\phi}}|^p|\sF_0)\leq C_{p,T,K}\left(\|\phi\|_{L^2_\mu}^p+|\phi(X_0)|^p\right).
\end{align}
Consequently, by \eqref{Pf2(ND)Le} and \eqref{Pf3(ND)Le} we obtain that for any $t\in[0,T]$,
\beg{align*}
\E(\zeta_t|\sF_0)\le K\int_0^t(\E(\zeta_s|\sF_0)+\E\zeta_s)\d s+C_{T,K,\tilde{K}}\chi(\ve,\phi),
\end{align*}
where
\beg{align*}
\chi(\ve,\phi)&=\ve^{\beta_0}\|\phi\|_{L^2_\mu}^{\beta_0+1}+\ve^{\gamma_0}\|\phi\|_{L^2_\mu}^{\gamma_0+1}+\ve\|\phi\|_{L^2_\mu}^2
+\ve^{\beta_0}\|\phi\|_{L^2_\mu}|\phi(X_0)|^{\beta_0}\cr
&\quad+\left[\ve^{\beta_0}\left(\|\phi\|_{L^2_\mu}^{\beta_0}+|\phi(X_0)|^{\beta_0}\right)+\ve\|\phi\|_{L^2_\mu}\right]|\phi(X_0)|.
\end{align*}
Taking the expectation on both sides and applying the Gronwall lemma, we obtain
\beg{align}\label{Pf4(ND)Le}
\E\zeta_t\le C_{T,K,\tilde{K}}\E\chi(\ve,\phi)\le C_{T,K,\tilde{K}}\ell(\ve,\phi)\|\phi\|_{L^2_\mu},
\end{align}
where $\ell(\ve,\phi)$ is given in \eqref{AD0Le2}.
Hence, this leads to our first claim \eqref{(ND)Le1}.

Next, we focus on proving \eqref{(ND)Le2}.
Applying the chain rule to $\zeta_t^2$ and using (A2) yield that for any $t\in[0,T]$,
\beg{align*}
&\d\zeta_t^2
\le 2K\zeta_t^2+2K\zeta_t\E\zeta_t\cr
&+3\tilde{K}(|\na_{\phi(X_0)}X_t^{\mu_{\ve,\phi}}|^2+|\na_{\phi(X_0)}X_t^\mu|^2)
\left(|X_t^{\mu_{\ve,\phi}}-X_t^\mu|^{\beta_0}+\W_2(\sL_{X_t^{\mu_{\ve,\phi}}},\sL_{X_t^\mu})\right)\cr
&+2\tilde{K}\zeta_t(\E|\na_{\phi(X_0)}X_t^\mu|^2)^\ff 1 2
\left(|X_t^{\mu_{\ve,\phi}}-X_t^\mu|^{\beta_0}+\W_2(\sL_{X_t^{\mu_{\ve,\phi}}},\sL_{X_t^\mu})+(\E|X_t^{\mu_{\ve,\phi}}-X_t^\mu|^{2\gamma_0})^\ff 1 2\right).
\end{align*}
Then, by the H\"{o}lder inequality we deduce that for any $t\in[0,T]$,
\beg{align*}
&\E(\zeta_t^2|\sF_0)\cr
&\le 2K\int_0^t\left[\E(\zeta_s^2|\sF_0)+\E(\zeta_s|\sF_0)\E\zeta_s\right]\d s\cr
&\quad+3\tilde{K}\int_0^t\left[\left(\E(|\na_{\phi(X_0)}X_s^{\mu_{\ve,\phi}}|^4|\sF_0)\right)^\ff 1 2+(\E\left(|\na_{\phi(X_0)}X_s^\mu|^4|\sF_0)\right)^\ff 1 2\right]\cr
&\qquad\quad\quad\times\left[\left(\E(|X_s^{\mu_{\ve,\phi}}-X_s^\mu|^{2\beta_0}\big|\sF_0)\right)^\ff 1 2+\W_2(\sL_{X_s^{\mu_{\ve,\phi}}},\sL_{X_s^\mu})\right]\d s\cr
&\quad+2\tilde{K}\int_0^t\E(\zeta_s^2|\sF_0)^\ff 1 2(\E|\na_{\phi(X_0)}X_s^\mu|^2)^\ff 1 2\cr
&\qquad\quad\quad\times\left[(\E(|X_s^{\mu_{\ve,\phi}}-X_s^\mu|^{2\beta_0}|\sF_0))^\ff 1 2+\W_2(\sL_{X_s^{\mu_{\ve,\phi}}},\sL_{X_s^\mu})+(\E|X_s^{\mu_{\ve,\phi}}-X_s^\mu|^{2\gamma_0})^\ff 1 2\right]\d s\cr
&\le \int_0^t\left[(2K+K^2+\tilde{K}^2)\E(\zeta_s^2|\sF_0)+(\E\zeta_s)^2\right]\d s\cr
&\quad+3\tilde{K}\int_0^t\left[\left(\E(|\na_{\phi(X_0)}X_s^{\mu_{\ve,\phi}}|^4|\sF_0)\right)^\ff 1 2+(\E\left(|\na_{\phi(X_0)}X_s^\mu|^4|\sF_0)\right)^\ff 1 2\right]\cr
&\qquad\quad\quad\times\left[\left(\E(|X_s^{\mu_{\ve,\phi}}-X_s^\mu|^{2\beta_0}\big|\sF_0)\right)^\ff 1 2
+(\E|X_s^{\mu_{\ve,\phi}}-X_s^\mu|^2)^\ff 1 2\right]\d s\cr
&\quad+3\int_0^t
\E|\na_{\phi(X_0)}X_s^\mu|^2
\left[\E(|X_s^{\mu_{\ve,\phi}}-X_s^\mu|^{2\beta_0}|\sF_0)+\E|X_s^{\mu_{\ve,\phi}}-X_s^\mu|^2+\E|X_s^{\mu_{\ve,\phi}}-X_s^\mu|^{2\gamma_0}\right]\d s.
\end{align*}
Combining this with \eqref{Pf2(ND)Le}, \eqref{Pf3(ND)Le} and  \eqref{Pf4(ND)Le} and applying the Gronwall lemma, we conclude that that for any $t\in[0,T]$,
\beg{align*}
\E(\zeta_t^2|\sF_0)\le C_{T,K,\tilde{K}}\left(\ell^2(\ve,\phi)\|\phi\|_{L^2_\mu}^2+\tilde{\chi}(\ve,\phi)\right),
\end{align*}
where
\beg{align*}
\tilde{\chi}(\ve,\phi)&=\left(\|\phi\|^2_{L^2_\mu}+|\phi(X_0)|^2\right)\left[\ve^{\beta_0}\left(\|\phi\|_{L^2_\mu}^{\beta_0}+|\phi(X_0)|^{\beta_0}\right)+\ve\|\phi\|_{L^2_\mu}\right]\cr
&\quad+\|\phi\|^2_{L^2_\mu}\left[\ve^{2\beta_0}\left(\|\phi\|_{L^2_\mu}^{2\beta_0}+|\phi(X_0)|^{2\beta_0}\right)+\ve^{2\gamma_0}\|\phi\|_{L^2_\mu}^{2\gamma_0}+\ve^2\|\phi\|^2_{L^2_\mu}\right].
\end{align*}
Letting $\tilde{\ell}_i(\ve,\phi), i=1,2,3$ are given respectively in \eqref{AD1Le2}-\eqref{AD3Le2},
we obtain the other assertion  \eqref{(ND)Le2}.
\end{proof}

Now, we are in the position to prove Theorem \ref{(ND)Th}.

\emph{Proof of Theorem \ref{(ND)Th}.}
For any  $\ve\in[0,1]$ and $s_0\in[0,T)$, let
\beg{align*}
\tilde{h}_{s_0}^{\ve,\phi}(t)
&=\int_{t\wedge s_0}^t\si^{-1}(s)\left[\ff 1{T-s_0} \na_{\phi(X_0)} X_s^{\mu_{\varepsilon,\phi}}\right.\cr
&\quad\left.+\ff {s-s_0}{T-s_0}\left(\E\<D^Lb(s,y,\cdot)(\sL_{X_s^{\mu_{\varepsilon,\phi}}})(X_s^{\mu_{\varepsilon,\phi}}), \na_{\phi(X_0)} X_s^{\mu_{\varepsilon,\phi}}\>\right)|_{y=X_s^{\mu_{\varepsilon,\phi}}}\right]\d s\cr
&=:\int_0^t\si^{-1}(s)\varrho_{\ve,s_0}(s)\mathrm{I}_{\{s>s_0\}}\d s,\ \ t\in[0,T].
\end{align*}
Owing to (A1) and (A2), one can verify that $\tilde{h}_{s_0}^{\ve,\phi}\in I_{0+}^{H+\ff 1 2}(L^2([0,T],\R^d))$,
which means that there exists $h_{s_0}^{\ve,\phi}\in\cH$ such that $R_H h_{s_0}^{\ve,\phi}=\tilde{h}_{s_0}^{\ve,\phi}$.
It is easy to see that $(R_H h_{s_0}^{\ve,\phi})(t)=0$ for all $t\in[0,s_0]$.
Moreover, applying the chain rule to $\ff {t-s_0}{T-s_0}\na_{\phi(X_0)} X_t^{\mu_{\varepsilon,\phi}}$ yields $\D_{R_H h_{s_0}^{\ve,\phi}}X_T^{\mu_{\varepsilon,\phi}}=\na_{\phi(X_0)} X_T^{\mu_{\varepsilon,\phi}}$.

Next, we intend to show that $K_H^{-1}(R_H h_{s_0}^{\ve,\phi})=K_H^*h_{s_0}^{\ve,\phi}\in L^2_a([0,T]\times\Omega,\R^d)$.
Then it follows from Remark \ref{Tra-Rem1} and Proposition \ref{Tra-Prp2} that $h_{s_0}^{\ve,\phi}\in\mathrm{Dom}\delta$ and
$\delta(h_{s_0}^{\ve,\phi})=\delta_W(K_H^*h_{s_0}^{\ve,\phi})=\int_0^T\<K_H^{-1}(R_H h_{s_0}^{\ve,\phi})(t),\d W_t\>$.

It is clear that the operator $K^{-1}_H$ preserves the adaptability property.
With the help of \eqref{InOp} and \eqref{FrDe}, we have
\beg{align}\label{addfor-de}
&K_H^{-1}\left(\int_0^\cdot\si^{-1}(s)\varrho_{\ve,s_0}(s)\mathrm{I}_{\{s>s_0\}}\d s\right)(t)\cr
&=t^{H-\frac{1}{2}}D^{H-\frac{1}{2}}_{0+}
\left[\cdot^{\frac{1}{2}-H}\si^{-1}(\cdot)\varrho_{\ve,s_0}(\cdot)\mathrm{I}_{\{\cdot>s_0\}}\right](t)\cr
&=\frac{H-\ff 1 2}{\Gamma(\frac{3}{2}-H)}
\Bigg[\ff {t^{\frac{1}{2}-H}\si^{-1}(t)\varrho_{\ve,s_0}(t)\mathrm{I}_{\{t>s_0\}}}{H-\ff 1 2}+\si^{-1}(t)\varrho_{\ve,s_0}(t)
\int_0^t\ff{\mathrm{I}_{\{t>s_0\}}-\mathrm{I}_{\{s>s_0\}}}{(t-s)^{\frac{1}{2}+H}}\d s\cr
&\qquad\qquad\qquad+t^{H-\frac{1}{2}}\si^{-1}(t)\varrho_{\ve,s_0}(t)\int_0^t\frac{t^{\frac{1}{2}-H}-s^{\frac{1}{2}-H}}{(t-s)^{\frac{1}{2}+H}}
\mathrm{I}_{\{s>s_0\}}\d s\cr
&\qquad\qquad\qquad+t^{H-\frac{1}{2}}\varrho_{\ve,s_0}(t)\int_0^t\frac{\si^{-1}(t)-\si^{-1}(s)}{(t-s)^{\frac{1}{2}+H}}s^{\ff 1 2-H}\mathrm{I}_{\{s>s_0\}}\d s\cr
&\qquad\qquad\qquad+t^{H-\frac{1}{2}}\int_0^t\frac{\varrho_{\ve,s_0}(t)-\varrho_{\ve,s_0}(s)}{(t-s)^{\frac{1}{2}+H}}\si^{-1}(s)s^{\ff 1 2-H}\mathrm{I}_{\{s>s_0\}}\d s\Bigg]\cr
&=:\frac{H-\ff 1 2}{\Gamma(\frac{3}{2}-H)}[I_1(t)+I_2(t)+I_3(t)+I_4(t)+I_5(t)].
\end{align}
From \eqref{Pf3(ND)Le}, it follows that
\beg{align*}
\sup_{s\in[0,T],\ve\in[0,1]}\E|\varrho_{\ve,s_0}(s)|^2\leq C_{s_0,T,K}\|\phi\|_{L^2_\mu}^2.
\end{align*}
Additionally, we have
\beg{align*}
\int_0^t\ff{\mathrm{I}_{\{t>s_0\}}-\mathrm{I}_{\{s>s_0\}}}{(t-s)^{\frac{1}{2}+H}}\d s
=\ff 1{H-\ff 1 2}\left((t-s_0)^{\ff 1 2-H}-t^{\ff 1 2-H}\right)\mathrm{I}_{\{t>s_0\}}
\end{align*}
and
\beg{align}\label{Pf6-(ND)Th}
\int_0^t\frac{s^{\frac{1}{2}-H}-t^{\frac{1}{2}-H}}{(t-s)^{\frac{1}{2}+H}}\d s
=t^{1-2H}\int_0^1\frac{r^{\frac{1}{2}-H}-1}{(1-r)^{\frac{1}{2}+H}}\d r<\infty.
\end{align}
These, together with (A2)(ii), imply that
\beg{align}\label{Pf2-(ND)Th}
&\E|I_1(t)|^2+\E|I_3(t)|^2\leq C_{s_0,T,K,H}t^{1-2H}\|\phi\|_{L^2_\mu}^2,\\
&\E|I_2(t)|^2\leq C_{s_0,T,K,H}(t-s_0)^{1-2H}\|\phi\|_{L^2_\mu}^2,\label{Pf3-(ND)Th}\\
&\E|I_4(t)|^2\leq C_{s_0,T,K,\tilde{K},H}t^{2\delta_0-2H+1}\|\phi\|_{L^2_\mu}^2,\label{Pf4-(ND)Th}
\end{align}
which means that $I_i\in L^2([0,T]\times\Omega,\R^d), i=1,\cdots,4.$\\
Before handing $I_5$, we set for any $\ve\in[0,1],t\in[0,T]$ and $y\in\R^d$,
\beg{align*}
\bar{b}^\ve(t,y):=D^Lb(t,y,\cdot)(\sL_{X_t^{\mu_{\ve,\phi}}})(X_t^{\mu_{\ve,\phi}}).
\end{align*}
By a direct calculation, we can reduce the integrability of $I_5$ to that of the following three terms in $L^2([0,T]\times\Omega,\R^d)$:
\beg{align*}
&t^{H-\frac{1}{2}}\int_0^t\frac{|\na_{\phi(X_0)} X_t^{\mu_{\varepsilon,\phi}}-\na_{\phi(X_0)} X_s^{\mu_{\varepsilon,\phi}}|}{(t-s)^{\frac{1}{2}+H}}s^{\ff 1 2-H}\d s,\cr
&t^{H-\frac{1}{2}}\int_0^t\frac{\left|\left(\E\<\bar{b}^\ve(t,y)-\bar{b}^\ve(s,z),
\na_{\phi(X_0)} X_t^{\mu_{\varepsilon,\phi}}\>\right)|_{y=X_t^{\mu_{\varepsilon,\phi}},z=X_s^{\mu_{\varepsilon,\phi}}}\right|}{(t-s)^{\frac{1}{2}+H}}s^{\ff 1 2-H}\d s,\cr
&t^{H-\frac{1}{2}}\int_0^t\frac{\left|\left(\E\<\bar{b}^\ve(s,z),
\na_{\phi(X_0)} X_t^{\mu_{\varepsilon,\phi}}-\na_{\phi(X_0)} X_s^{\mu_{\varepsilon,\phi}}\>\right)|_{z=X_s^{\mu_{\varepsilon,\phi}}}\right|}
{(t-s)^{\frac{1}{2}+H}}s^{\ff 1 2-H}\d s.
\end{align*}
Along the same lines as in Step 1 of Theorem \ref{(EX)The1}, we have for any $p\in(1/H,2]$,
\beg{align}\label{Pf20-(ND)Th}
\sup_{\ve\in[0,1]}\E|X_t^{\mu_{\varepsilon,\phi}}-X_s^{\mu_{\varepsilon,\phi}}|^p\leq C_{p,T,K,H}|t-s|^{pH},\ \  t,s\in[0,T].
\end{align}
Furthermore, by \eqref{Pf3(ND)Le}, it is easy to see that for any $p>0$ and $s,t\in[0,T]$,
\beg{align}\label{Pf10-(ND)Th}
\E(|\na_{\phi(X_0)} X_t^{\mu_{\varepsilon,\phi}}-\na_{\phi(X_0)} X_s^{\mu_{\varepsilon,\phi}}|^p|\sF_0)
\leq C_{p,T,K}\left(\|\phi\|_{L^2_\mu}^p+|\phi(X_0)|^p\right)|t-s|^p.
\end{align}
Then, combining these with (A2) implies
\beg{align}\label{Pf5'(ND)Th}
\E|I_5(t)|^2&\le C_{s_0,T,K,\ti K,H}\Bigg[t^{3-2H}+t^{2\alpha_0-2H+1}+t^{2\gamma_0H-2H+1}+t\cr
&\qquad\qquad\qquad+t^{2H-1}\E\bigg(\int_0^t\frac{|X_t^{\mu_{\varepsilon,\phi}}-X_s^{\mu_{\varepsilon,\phi}}|^{\beta_0}}{(t-s)^{\frac{1}{2}+H}}s^{\ff 1 2-H}\d s\bigg)^2\Bigg]\|\phi\|_{L^2_\mu}^2.
\end{align}
Note that there hold
\beg{align*}
&\sup_{r\in[0,T]}|X_r^{\mu_{\ve,\phi}}|\le C_{T,K,H}\left(1+\|\mathrm{Id}+\ve\phi\|_{L^2_\mu}+|X_0+\ve\phi(X_0)|+\left\|\int_0^\cdot\si(r)\d B^H_r\right\|_\infty\right)
\end{align*}
and
\beg{align*}
&\E\bigg(\sup_{r\in[0,T]}|X_r^{\mu_{\ve,\phi}}|^2\bigg)\le C_{T,K,H}\left(1+\|\mathrm{Id}+\ve\phi\|_{L^2_\mu}^2\right),
\end{align*}
where $\|\int_0^\cdot\si(r)\d B^H_r\|_\infty:=\sup_{t\in[0,T]}|\int_0^t\si(r)\d B^H_r|$.\\
Then, it follows from (A1) that
\beg{align}\label{term4-4}
&\left|\int_s^tb(r,X_r^{\mu_{\ve,\phi}},\sL_{X_r^{\mu_{\ve,\phi}}})\d r\right|\cr
&\le K\left(1+\sup_{r\in[0,T]}|X_r^{\mu_{\ve,\phi}}|+\Big(\E\sup_{r\in[0,T]}|X_r^{\mu_{\ve,\phi}}|^2\Big)^\ff1 2\right)(t-s)\cr
&\le C_{T,K,H}\left(1+\|\mathrm{Id}+\ve\phi\|_{L^2_\mu}+|X_0+\ve\phi(X_0)|+\left\|\int_0^\cdot\si(r)\d B^H_r\right\|_\infty\right)(t-s).
\end{align}
Consequently, this implies
\beg{align}\label{Add1-I5}
&t^{2H-1}\E\bigg(\int_0^t\frac{|X_t^{\mu_{\varepsilon,\phi}}-X_s^{\mu_{\varepsilon,\phi}}|^{\beta_0}}{(t-s)^{\frac{1}{2}+H}}s^{\ff 1 2-H}\d s\bigg)^2\cr
&\le2t^{2H-1}\E\bigg(\int_0^t\ff {|\int_s^tb(r,X_r^{\mu_{\ve,\phi}},\sL_{X_r^{\mu_{\ve,\phi}}})\d r|^{\be_0}}{(t-s)^{\ff 1 2+H}}s^{\ff 1 2-H}\d s\bigg)^2\cr
&\quad+2t^{2H-1}\E\bigg(\int_0^t\ff {|\int_s^t\si(r)\d B^H_r|^{\be_0}}{(t-s)^{\ff 1 2+H}}s^{\ff 1 2-H}\d s\bigg)^2\nonumber\\
&\le C_{T,K,H}\Bigg(1+\|\mathrm{Id}+\ve\phi\|_{L^2_\mu}^{2\be_0}+\E\left\|\int_0^\cdot\si(r)\d B^H_r\right\|_\infty^{2\be_0}\Bigg)t^{1+2(\be_0-H)}\nonumber\\
&\quad+C_H\E\left(\left\|\int_0^\cdot\si(r)\d B^H_r\d r\right\|^{2\be_0}_{H-\varsigma_0}\right)t^{1+2(H-\varsigma_0)\be_0-2H},
\end{align}
where we use the H\"{o}lder continuity of $\int_0^\cdot\si(r)\d B^H_r$ of order $H-\varsigma_0$ with $\varsigma_0\in(0,1/2)$ and
\beg{align*}
\Big\|\int_0^\cdot\si(r)\d B^H_r\Big\|_{H-\varsigma_0}:=\sup_{0\leq s<t\leq T}\ff {\left|\int_0^t\si(r)\d B^H_r-\int_0^s\si(r)\d B^H_r\right|}{|t-s|^{H-\varsigma_0}}.
\end{align*}
Plugging \eqref{Add1-I5} into \eqref{Pf5'(ND)Th} yields that $I_5\in L^2([0,T]\times\Omega,\R^d)$.
Then we get the desired claim.

Since
\beg{align*}
\E\delta^2(h_{s_0}^{\ve,\phi})=\E\delta_W^2(K_H^*h_{s_0}^{\ve,\phi})=\int_0^T\E|K_H^{-1}(R_H h_{s_0}^{\ve,\phi})(t)|^2\d t,
\end{align*}
by \eqref{Pf2-(ND)Th}, \eqref{Pf3-(ND)Th}, \eqref{Pf4-(ND)Th}, \eqref{Pf5'(ND)Th} and \eqref{Add1-I5}, one has that $\int_0^1\left(\E\delta^2(h_{s_0}^{\tau,\phi})\right)^\ff 1 2\d\tau<\infty$ and
$\E\delta^2(h^\phi)\leq C_{T,K,\tilde{K},H}\|\phi\|^2_{L^2_\mu}$.

Finally, we shall estimate $\E|\delta(h^{\ve,\phi})-\delta(h^\phi)|$.
As before, we write $\varrho_{\ve}=\varrho_{\ve,0}$ and $\varrho=\varrho_{0,0}$ for simplicity.
Using the linearity of the operator $K_H^{-1}$ and applying the B.D.G. inequality and the H\"{o}lder inequality, we have
\beg{align}\label{Pf18-(ND)Th}
\E|\delta(h^{\ve,\phi})-\delta(h^\phi)|&=\E\left|\int_0^T\<K_H^{-1}(R_H h^{\ve,\phi})(t),\d W_t\>-\int_0^T\<K_H^{-1}(R_H h^\phi)(t),\d W_t\>\right|\nonumber\\
&\le\E\left(\int_0^T|K_H^{-1}(R_H h^{\ve,\phi}-R_H h^{\phi})(t)|^2\d t\right)^\ff 1 2 \nonumber\\
&=\E\left(\E\left(\left(\int_0^T|K_H^{-1}(R_H h^{\ve,\phi}-R_H h^{\phi})(t)|^2\d t\right)^\ff 1 2\big|\sF_0\right)\right) \nonumber\\
&\le\E\left(\E\left(\int_0^T|K_H^{-1}(R_H h^{\ve,\phi}-R_H h^{\phi})(t)|^2\d t\big|\sF_0\right)\right)^\ff 1 2 \nonumber\\
&=\E\left(\int_0^T\E\left(\left|K_H^{-1}\left(\int_0^\cdot\si^{-1}(s)(\varrho_{\ve}-\varrho)(s)\d s\right)(t)\right|^2\big|\sF_0\right)\d t\right)^\ff 1 2.
\end{align}
By \eqref{InOp} and \eqref{FrDe} again, we have
\beg{align}\label{ForRem}
&K_H^{-1}\left(\int_0^\cdot\si^{-1}(s)(\varrho_{\ve}-\varrho)(s)\d s\right)(t)=t^{H-\frac{1}{2}}D^{H-\frac{1}{2}}_{0+}
\left[\cdot^{\frac{1}{2}-H}\si^{-1}(\cdot)(\varrho_{\ve}-\varrho)(\cdot)\right](t)\cr
&=\frac{H-\ff 1 2}{\Gamma(\frac{3}{2}-H)}
\Bigg[\ff {t^{\frac{1}{2}-H}\si^{-1}(t)(\varrho_{\ve}-\varrho)(t)}{H-\ff 1 2}
+t^{H-\frac{1}{2}}\si^{-1}(t)(\varrho_{\ve}-\varrho)(t)\int_0^t\frac{t^{\frac{1}{2}-H}-s^{\frac{1}{2}-H}}{(t-s)^{\frac{1}{2}+H}}
\d s\cr
&\qquad\qquad\qquad+t^{H-\frac{1}{2}}(\varrho_{\ve}-\varrho)(t)\int_0^t\frac{\si^{-1}(t)-\si^{-1}(s)}{(t-s)^{\frac{1}{2}+H}}s^{\ff 1 2-H}\d s\cr
&\qquad\qquad\qquad+t^{H-\frac{1}{2}}\int_0^t\frac{(\varrho_{\ve}-\varrho)(t)-(\varrho_{\ve}-\varrho)(s)}{(t-s)^{\frac{1}{2}+H}}\si^{-1}(s)s^{\ff 1 2-H}\d s\Bigg]\cr
&=:\frac{H-\ff 1 2}{\Gamma(\frac{3}{2}-H)}[J_1(t)+J_2(t)+J_3(t)+J_4(t)].
\end{align}
Owing to (A2)(ii) and \eqref{Pf6-(ND)Th}, we get
\beg{align*}
\sum_{i=1}^3\E(|J_i(t)|^2|\sF_0)\le C_{T,\tilde{K},H}\left(t^{1-2H}+t^{2\delta_0-2H+1}\right)\E\left(|(\varrho_{\ve}-\varrho)(t)|^2|\sF_0\right),
\end{align*}
which leads to
\beg{align}\label{Pf7-(ND)Th}
\E\left(\int_0^T\sum_{i=1}^3\E(|J_i(t)|^2|\sF_0)\d t\right)^\ff 1 2\le C_{T,\tilde{K},H}\E\left(\sup_{t\in[0,T]}\E\left(|(\varrho_{\ve}-\varrho)(t)|^2|\sF_0\right)\right)^\ff 1 2.
\end{align}
Note that by (A1) and (A2), we obtain for any $t\in[0,T]$,
\beg{align*}
\E\left(|(\varrho_{\ve}-\varrho)(t)|^2|\sF_0\right)
&\le\ff 3{T^2}\E(|\na_{\phi(X_0)}X_t^{\mu_{\ve,\phi}}-\na_{\phi(X_0)}X_t^\mu|^2|\sF_0)\cr
&\quad+6\tilde{K}^2\left((\E(|X_t^{\mu_{\ve,\phi}}-X_t^\mu|^{\beta_0}|\sF_0)+\W_2(\sL_{X_t^{\mu_{\ve,\phi}}},\sL_{X_t^\mu})\right)^2(\E|\na_{\phi(X_0)}X_t^\mu|)^2\cr
&\quad+6\tilde{K}^2\left(\E(|X_t^{\mu_{\ve,\phi}}-X_t^\mu|^{\gamma_0}|\na_{\phi(X_0)}X_t^\mu|)\right)^2\cr
&\quad+3K^2(\E|\na_{\phi(X_0)}X_t^{\mu_{\ve,\phi}}-\na_{\phi(X_0)}X_t^\mu|)^2\cr
&\le C_{T,K,\tilde{K}}\Big[\left(\ell^2(\ve,\phi)+\tilde{\ell}_1^2(\ve,\phi)+\ve^{2\beta_0}\|\phi\|_{L^2_\mu}^{2\beta_0}+\ve^2\|\phi\|_{L^2_\mu}^2+\ve^{2\gamma_0}\|\phi\|_{L^2_\mu}^{2\gamma_0}\right)\|\phi\|_{L^2_\mu}^2\cr
&\qquad\qquad\quad+(\tilde{\ell}_2^2(\ve,\phi)+\ve^{2\beta_0}|\phi(X_0)|^{2\beta_0})\|\phi\|_{L^2_\mu}^2+\tilde{\ell}_3^2(\ve,\phi)|\phi(X_0)|^2\Big],
\end{align*}
where the last inequality is due to \eqref{Pf2(ND)Le}, \eqref{Pf3(ND)Le} and Lemma \ref{(ND)Le}.
Then, combining this with Remark \ref{(ND)Re1} and \eqref{Pf7-(ND)Th} yields
\beg{align}\label{Pf8-(ND)Th}
\lim_{\ve\ra0}\E\left(\int_0^T\sum_{i=1}^3\E(|J_i(t)|^2|\sF_0)\d t\right)^\ff 1 2=0
\end{align}
and
\beg{align}\label{Pf9-(ND)Th}
\lim_{\|\phi\|_{L^2_\mu}\ra 0}\sup\limits_{\ve\in(0,1]}\ff{\E\left(\int_0^T\sum_{i=1}^3\E(|J_i(t)|^2|\sF_0)\d t\right)^\ff 1 2}{\|\phi\|_{L^2_\mu}}=0.
\end{align}
For the term associated with $J_4(t)$, observe first that for any $0\le s<t\le T$,
\beg{align}\label{Pf17-(ND)Th}
(\varrho_{\ve}-\varrho)(t)-(\varrho_{\ve}-\varrho)(s)=\sum_{i=1}^6\Theta_i(t,s),
\end{align}
where
\beg{align*}
&\Theta_1(t,s)=\ff 1 T\left[\left(\na_{\phi(X_0)}X_t^{\mu_{\ve,\phi}}-\na_{\phi(X_0)}X_s^{\mu_{\ve,\phi}}\right)
-\left(\na_{\phi(X_0)}X_t^\mu-\na_{\phi(X_0)}X_s^\mu\right)\right],\cr
&\Theta_2(t,s)=\ff {t-s}T\left[\left(\E\<\bar{b}^\ve(t,y), \na_{\phi(X_0)} X_t^{\mu_{\ve,\phi}}\>\right)|_{y=X_t^{\mu_{\ve,\phi}}}
-\left(\E\<\bar{b}(t,\ti y), \na_{\phi(X_0)} X_t^\mu\>\right)|_{\ti y=X_t^\mu}\right],\cr
&\Theta_3(t,s)=\ff s T\left(\E\left\<\bar{b}^\ve(t,y),\left(\na_{\phi(X_0)}X_t^{\mu_{\ve,\phi}}-\na_{\phi(X_0)}X_s^{\mu_{\ve,\phi}}\right)
-\left(\na_{\phi(X_0)}X_t^\mu-\na_{\phi(X_0)}X_s^\mu\right)\right\>\right)|_{y=X_t^{\mu_{\ve,\phi}}},\cr
&\Theta_4(t,s)=\ff s T\left(\E\<\bar{b}^\ve(t,y)-\bar{b}^\ve(s,z),\na_{\phi(X_0)}X_s^{\mu_{\ve,\phi}}-\na_{\phi(X_0)} X_s^\mu\>\right)|_{y=X_t^{\mu_{\ve,\phi}},z=X_s^{\mu_{\ve,\phi}}},\cr
&\Theta_5(t,s)=\ff s T\left(\E\<\bar{b}^\ve(s,z)-\bar{b}(s,\ti z),\na_{\phi(X_0)} X_t^\mu-\na_{\phi(X_0)} X_s^\mu\>\right)|_{z=X_s^{\mu_{\ve,\phi}},\ti z=X_s^\mu},\cr
&\Theta_6(t,s)=\ff s T\left(\E\<(\bar{b}^\ve(t,y)-\bar{b}^\ve(s,z))-(\bar{b}(t,\ti y)-\bar{b}(s,\ti z)),\na_{\phi(X_0)} X_t^\mu\>\right)|_{y=X_t^{\mu_{\ve,\phi}},z=X_s^{\mu_{\ve,\phi}},\ti y=X_t^\mu,\ti z=X_s^\mu}.
\end{align*}
Owing to  (A1), (A2), \eqref{Pf2(ND)Le}, \eqref{Pf3(ND)Le}, \eqref{Pf10-(ND)Th} and Lemma \ref{(ND)Le}, one gets that
\beg{align}\label{Pf11-(ND)Th}
&\E\left(t^{2H-1}\left(\int_0^t\ff {|\Theta_1(t,s)|s^{\ff 1 2-H}}{(t-s)^{\ff 1 2+H}}\d s\right)^2\big|\sF_0\right)\nonumber\\
&\le  C_{T,K,\tilde{K},H}t^{3-2H}
\left[\tilde{\ell}_1^2(\ve,\phi)\|\phi\|_{L^2_\mu}^2+\tilde{\ell}_2^2(\ve,\phi)\|\phi\|_{L^2_\mu}^2+(\tilde{\ell}_3^2(\ve,\phi)+\tilde{\ell}_3^4(\ve,\phi))|\phi(X_0)|^2\right],\\
&\E\left(t^{2H-1}\left(\int_0^t\ff {|\Theta_2(t,s)|s^{\ff 1 2-H}}{(t-s)^{\ff 1 2+H}}\d s\right)^2\big|\sF_0\right)\label{Pf12-(ND)Th}\nonumber\\
&\le  C_{T,K,\tilde{K},H}t^{3-2H}
\left(\ell^2(\ve,\phi)\|\phi\|_{L^2_\mu}^2+\ve^{2\beta_0}|\phi(X_0)|^{2\beta_0}\|\phi\|_{L^2_\mu}^2\right),\\
&\E\left(t^{2H-1}\left(\int_0^t\ff {|\Theta_3(t,s)|s^{\ff 1 2-H}}{(t-s)^{\ff 1 2+H}}\d s\right)^2\big|\sF_0\right)\label{Pf13-(ND)Th}\\
&\le  C_{T,K,\tilde{K},H}t^{3-2H}\ell^2(\ve,\phi)\|\phi\|_{L^2_\mu}^2,\nonumber\\
&\E\left(t^{2H-1}\left(\int_0^t\ff {|\Theta_5(t,s)|s^{\ff 1 2-H}}{(t-s)^{\ff 1 2+H}}\d s\right)^2\big|\sF_0\right)\label{Pf14-(ND)Th}\nonumber\\
&\le  C_{T,K,\tilde{K},H}t^{3-2H}\left[\ve^2\|\phi\|_{L^2_\mu}^2+\ve^{2\ga_0}\|\phi\|_{L^2_\mu}^{2\ga_0}+
\ve^{2\be_0}\left(\|\phi\|_{L^2_\mu}^{2\be_0}+|\phi(X_0)|^{2\be_0}\right)\right]
\|\phi\|_{L^2_\mu}^2.\nonumber\\
\end{align}
For $\Theta_4(t,s)$, by (A2)(i), \eqref{(ND)Le1} and  \eqref{Pf20-(ND)Th} we first have
\beg{align}\label{term4-2}
|\Theta_4(t,s)|&\le C_{T,K,\tilde{K},H}\left[(t-s)^{\al_0}+(t-s)^H+|X_t^{\mu_{\ve,\phi}}-X_s^{\mu_{\ve,\phi}}|^{\be_0}\right]\ell(\ve,\phi)\|\phi\|_{L^2_\mu}\cr
&\quad+\tilde{K}\E\left(|X_t^{\mu_{\ve,\phi}}-X_s^{\mu_{\ve,\phi}}|^{\ga_0}\cdot|\na_{\phi(X_0)}X_s^{\mu_{\ve,\phi}}-\na_{\phi(X_0)} X_s^\mu|\right)\cr
&\le C_{T,K,\tilde{K},H}\left[(t-s)^{\al_0}+(t-s)^H+|X_t^{\mu_{\ve,\phi}}-X_s^{\mu_{\ve,\phi}}|^{\be_0}\right]\ell(\ve,\phi)\|\phi\|_{L^2_\mu}\cr
&\quad+\tilde{K}\E\left(\left|\int_s^tb(r,X_r^{\mu_{\ve,\phi}},\sL_{X_r^{\mu_{\ve,\phi}}})\d r\right|^{\ga_0}\cdot|\na_{\phi(X_0)}X_s^{\mu_{\ve,\phi}}-\na_{\phi(X_0)} X_s^\mu|\right)\cr
&\quad+\tilde{K}\E\left(\left|\int_s^t\si(r)\d B^H_r\right|^{\ga_0}\cdot|\na_{\phi(X_0)}X_s^{\mu_{\ve,\phi}}-\na_{\phi(X_0)} X_s^\mu|\right).
\end{align}
Next, we focus on dealing with the last two terms of the right-hand side of \eqref{term4-2}.
Using \eqref{term4-4}, \eqref{(ND)Le1}, \eqref{(ND)Le2} and the fact that $B^H$ is independent of $\sF_0$, we obtain
\beg{align}\label{term4-1}
&\E\left(\left|\int_s^tb(r,X_r^{\mu_{\ve,\phi}},\sL_{X_r^{\mu_{\ve,\phi}}})\d r\right|^{\ga_0}\cdot|\na_{\phi(X_0)}X_s^{\mu_{\ve,\phi}}-\na_{\phi(X_0)} X_s^\mu|\right)\cr
&\le C_{T,K,\tilde{K},H}\Bigg\{\left(1+\|\mathrm{Id}+\ve\phi\|_{L^2_\mu}^{\ga_0}\right)\ell(\ve,\phi)\|\phi\|_{L^2_\mu}\cr
&\qquad\qquad\quad+\left(\E\left\|\int_0^\cdot\si(r)\d B^H_r\right\|_\infty^{2\ga_0}\right)^\ff 1 2\left[\tilde{\ell}_1(\ve,\phi)+\E\tilde{\ell}_2(\ve,\phi)+\left(\E\tilde{\ell}_3^2(\ve,\phi)\right)^\ff 1 2\right]\|\phi\|_{L^2_\mu}\cr
&\qquad\qquad\quad+\E\left(|X_0+\ve\phi(X_0)|^{\ga_0}\cdot|\na_{\phi(X_0)}X_s^{\mu_{\ve,\phi}}-\na_{\phi(X_0)} X_s^\mu|\right)\Bigg\}(t-s)^{\ga_0}.
\end{align}
Observe that by \eqref{(ND)Le2}, we derive
\beg{align}\label{term4-5}
&\E\left(|X_0+\ve\phi(X_0)|^{\ga_0}\cdot|\na_{\phi(X_0)}X_s^{\mu_{\ve,\phi}}-\na_{\phi(X_0)} X_s^\mu|\right)\cr
&\le C_{T,K,\tilde{K}}\E\left[|X_0+\ve\phi(X_0)|^{\ga_0}
\left(\tilde{\ell}_1(\ve,\phi)\|\phi\|_{L^2_\mu}+\tilde{\ell}_2(\ve,\phi)\|\phi\|_{L^2_\mu}+\tilde{\ell}_3(\ve,\phi)|\phi(X_0)|\right)\right]\cr
&\le C_{T,K,\tilde{K}}\bigg[\|\mathrm{Id}+\ve\phi\|_{L^2_\mu}^{\ga_0}\left(\tilde{\ell}_1(\ve,\phi)+\left(\E\tilde{\ell}^2_2(\ve,\phi)\right)^\ff 1 2+\ve^{\ff {\beta_0} 2}\|\phi\|_{L^2_\mu}^{\ff {\beta_0} 2}+\ve^{\ff 1 2}\|\phi\|^{\ff 1 2}_{L^2_\mu}\right)\|\phi\|_{L^2_\mu}\cr
&\qquad\qquad\quad+\ve^{\ff {\beta_0}2}\E\left(|X_0+\ve\phi(X_0)|^{\ga_0}\cdot|\phi(X_0)|^{1+\ff {\beta_0} 2}\right)\bigg]\cr
&\le C_{T,K,\tilde{K}}\|\mathrm{Id}+\ve\phi\|_{L^2_\mu}^{\ga_0}\left(\tilde{\ell}_1(\ve,\phi)+\left(\E\tilde{\ell}^2_2(\ve,\phi)\right)^\ff 1 2+\ve^{\ff {\beta_0} 2}\|\phi\|_{L^2_\mu}^{\ff {\beta_0} 2}+\ve^{\ff 1 2}\|\phi\|^{\ff 1 2}_{L^2_\mu}\right)\|\phi\|_{L^2_\mu}\cr
&= C_{T,K,\tilde{K}}\|\mathrm{Id}+\ve\phi\|_{L^2_\mu}^{\ga_0}\left(\tilde{\ell}_1(\ve,\phi)+\left(\E\tilde{\ell}^2_2(\ve,\phi)\right)^\ff 1 2\right)\|\phi\|_{L^2_\mu},
\end{align}
where we use the H\"{o}lder inequality with $\ff {2+\be_0} 4+\ff {2-\be_0} 4=1$ and the relation  $(1-\frac{1}{2H}\le)\ga_0\le 1-\ff {\be_0} 2$ in the last inequality.
Note that if $\ga_0\in(1-\ff {\be_0} 2,1]$,
we may choose $\ti\ga_0\in[1-\frac{1}{2H}, 1-\ff {\be_0} 2]$ to replace such $\ga_0$ in the first inequality of \eqref{term4-2} due to the boundedness of $D^Lb$.
In this case, \eqref{Pf15-(ND)Th} below holds with $\ga_0$ replaced by $\ti\ga_0$,
which also implies the desired convergence of the term involved $\Theta_4$.\\
Substituting \eqref{term4-5} into \eqref{term4-1} and recalling that $\ve\in[0,1]$  imply
\beg{align}\label{term4-3}
&\E\left(\left|\int_s^tb(r,X_r^{\mu_{\ve,\phi}},\sL_{X_r^{\mu_{\ve,\phi}}})\d r\right|^{\ga_0}\cdot|\na_{\phi(X_0)}X_s^{\mu_{\ve,\phi}}-\na_{\phi(X_0)} X_s^\mu|\right)\cr
&\le C_{T,K,\tilde{K},H}(t-s)^{\ga_0}\bigg(\ell(\ve,\phi)+\tilde{\ell}_1(\ve,\phi)+\left(\E\tilde{\ell}^2_2(\ve,\phi)\right)^\ff 1 2+\left(\E\tilde{\ell}_3^2(\ve,\phi)\right)^\ff 1 2\bigg)\|\phi\|_{L^2_\mu},
\end{align}
For the other term, applying the fact that $B^H$ is independent of $\sF_0$ again and \eqref{(ND)Le2}, one sees that
\beg{align*}
&\E\left(\left|\int_s^t\si(r)\d B^H_r\right|^{\ga_0}\cdot|\na_{\phi(X_0)}X_s^{\mu_{\ve,\phi}}-\na_{\phi(X_0)} X_s^\mu|\right)\cr
&\le\E\left\{\left[\E\left(\left|\int_s^t\si(r)\d B^H_r\right|^{2\ga_0}\big|\sF_0\right)\right]^\ff 1 2\cdot\left[\E\left(|\na_{\phi(X_0)}X_s^{\mu_{\ve,\phi}}-\na_{\phi(X_0)} X_s^\mu|^2|\sF_0\right)\right]^\ff 1 2\right\}\cr
&\le C_{T,K,\tilde{K},H}(t-s)^{\ga_0H}
\left(\tilde{\ell}_1(\ve,\phi)+\E\tilde{\ell}_2(\ve,\phi)+\left(\E\tilde{\ell}_3^2(\ve,\phi)\right)^\ff 1 2\right)\|\phi\|_{L^2_\mu}.
\end{align*}
Plugging this and \eqref{term4-3} into \eqref{term4-2}, we arrive at
\beg{align*}
|\Theta_4(t,s)|
&\le C_{T,K,\tilde{K},H}\bigg[|X_t^{\mu_{\ve,\phi}}-X_s^{\mu_{\ve,\phi}}|^{\be_0}\ell(\ve,\phi)\cr
&\quad+(t-s)^{\al_0\wedge(\ga_0H)}\left(\ell(\ve,\phi)+\tilde{\ell}_1(\ve,\phi)+\left(\E\tilde{\ell}^2_2(\ve,\phi)\right)^\ff 1 2+\left(\E\tilde{\ell}_3^2(\ve,\phi)\right)^\ff 1 2\right)\bigg]\|\phi\|_{L^2_\mu}.
\end{align*}
Hence, combining this with \eqref{Add1-I5} and the fact that $B^H$ is independent of $\sF_0$ again leads to
\beg{align}\label{Pf15-(ND)Th}
&\E\Bigg(t^{2H-1}\left(\int_0^t\ff {|\Theta_4(t,s)|s^{\ff 1 2-H}}{(t-s)^{\ff 1 2+H}}\d s\right)^2\big|\sF_0\Bigg)\nonumber\\
&\le C_{T,K,\tilde{K},H}\left(1+\|\mathrm{Id}+\ve\phi\|_{L^2_\mu}^{2\be_0}+\E\bigg\|\int_0^\cdot\si(r)\d B^H_r\bigg\|_\infty^{2\be_0}\right)t^{1+2(\be_0-H)}\ell^2(\ve,\phi)\|\phi\|^2_{L^2_\mu}\nonumber\\
&\quad+C_{T,K,\tilde{K},H}\E\left(\left\|\int_0^\cdot\si(r)\d B^H_r\d r\right\|^{2\be_0}_{H-\varsigma_0}\right)t^{1+2(H-\varsigma_0)\be_0-2H}\ell^2(\ve,\phi)\|\phi\|^2_{L^2_\mu}\nonumber\\
&\quad+C_{T,K,\tilde{K},H}
t^{1+2(\al_0\wedge(\ga_0H)-H)}\left(\ell^2(\ve,\phi)+\tilde{\ell}^2_1(\ve,\phi)+\E\tilde{\ell}^2_2(\ve,\phi)+\E\tilde{\ell}_3^2(\ve,\phi)\right)\|\phi\|^2_{L^2_\mu}\cr
&\le C_{T,K,\tilde{K},H}(t^{1+2(\be_0-H)}+t^{1+2(H-\varsigma_0)\be_0-2H})\ell^2(\ve,\phi)\|\phi\|^2_{L^2_\mu}\nonumber\\
&\quad+C_{T,K,\tilde{K},H}
t^{1+2(\al_0\wedge(\ga_0H)-H)}\left(\ell^2(\ve,\phi)+\tilde{\ell}^2_1(\ve,\phi)+\E\tilde{\ell}^2_2(\ve,\phi)+\E\tilde{\ell}_3^2(\ve,\phi)\right)\|\phi\|^2_{L^2_\mu}.
\end{align}
As far as $\Theta_6(t,s)$ is concerned,
using (A3) and Lemma \ref{FoLD}, we derive that for any $\ve\in[0,1],s,t\in[0,T]$ and $y,z\in\R^d$,
\beg{align*}
&\bar{b}^\ve(t,y)-\bar{b}^\ve(s,z)=D^Lb(t,y,\cdot)(\sL_{X_t^{\mu_{\ve,\phi}}})(X_t^{\mu_{\ve,\phi}})
-D^Lb(s,z,\cdot)(\sL_{X_s^{\mu_{\ve,\phi}}})(X_s^{\mu_{\ve,\phi}})\cr
&=\int_0^1\ff \d {\d\th}D^Lb(\th_{s,t},y,\cdot)(\sL_{X_t^{\mu_{\ve,\phi}}})(X_t^{\mu_{\ve,\phi}})\d\th\cr
&\quad+\int_0^1\ff \d {\d\th}D^Lb(s,z+\th(y-z),\cdot)(\sL_{X_t^{\mu_{\ve,\phi}}})(X_t^{\mu_{\ve,\phi}})\d\th\cr
&\quad+\int_0^1\ff \d {\d\th}D^Lb(s,z,\cdot)(\sL_{X^{\ve,\phi}_{s,t}(\th)})(X_t^{\mu_{\ve,\phi}})\d\th\cr
&\quad+\int_0^1\ff \d {\d\th}D^Lb(s,z,\cdot)(\sL_{X_s^{\mu_{\ve,\phi}}})(X^{\ve,\phi}_{s,t}(\th))\d\th\cr
&=\int_0^1\partial_{\th_{s,t}}(D^Lb(\cdot,y,\cdot)(\sL_{X_t^{\mu_{\ve,\phi}}})(X_t^{\mu_{\ve,\phi}}))(\th_{s,t})(t-s)\d\th\cr
&\quad+\int_0^1 \na(D^Lb(s,\cdot,\cdot)(\sL_{X_t^{\mu_{\ve,\phi}}})(X_t^{\mu_{\ve,\phi}}))(z+\th(y-z))(y-z)\d\th\cr
&\quad+\int_0^1\left(\E\<D^L(D^Lb(s,z,\cdot)(\cdot)(u))(\sL_{X^{\ve,\phi}_{s,t}(\th)})(X^{\ve,\phi}_{s,t}(\th)),
X_t^{\mu_{\ve,\phi}}-X_s^{\mu_{\ve,\phi}}\>\right)|_{u=X_t^{\mu_{\ve,\phi}}}
\d\th\cr
&\quad+\int_0^1 \na(D^Lb(s,z,\cdot)(\sL_{X_s^{\mu_{\ve,\phi}}})(\cdot))(X^{\ve,\phi}_{s,t}(\th))(X_t^{\mu_{\ve,\phi}}-X_s^{\mu_{\ve,\phi}})\d\th,
\end{align*}
where for any $\th\in[0,1]$, $\th_{s,t}:=s+\th(t-s)$ and $X^{\ve,\phi}_{s,t}(\th):=X_s^{\mu_{\ve,\phi}}+\th(X_t^{\mu_{\ve,\phi}}-X_s^{\mu_{\ve,\phi}})$.\\
Then by (A1), (A3) and \eqref{Pf20-(ND)Th}, we have
\beg{align*}
|\Theta_6(t,s)|\le C_{T,K,\bar K,H} \left[\sum_{i=1}^4\Lambda_i+|(X_t^{\mu_{\ve,\phi}}-X_s^{\mu_{\ve,\phi}})-(X_t^\mu-X_s^\mu)|+\ve(t-s)\|\phi\|_{L^2_\mu}\right]\|\phi\|_{L^2_\mu},
\end{align*}
where
\beg{align*}
\Lambda_1&:=(t-s)
\bigg(\E\int_0^1\bigg|\partial_{\th_{s,t}}(D^Lb(\cdot,y,\cdot)(\sL_{X_t^{\mu_{\ve,\phi}}})(X_t^{\mu_{\ve,\phi}}))(\th_{s,t})\\
&\qquad\qquad\qquad\quad-\partial_{\th_{s,t}}(D^Lb(\cdot,\ti y,\cdot)(\sL_{X_t^\mu})(X_t^\mu))(\th_{s,t})\bigg|^2\d\th\bigg)^\ff 1 2\big|_{y=X_t^{\mu_{\ve,\phi}},\ti y=X_t^\mu},\cr
\Lambda_2&:=\bigg(\E\int_0^1 |\na(D^Lb(s,\cdot,\cdot)(\sL_{X_t^{\mu_{\ve,\phi}}})(X_t^{\mu_{\ve,\phi}}))(z+\th(y-z))\\
&\qquad\qquad-\na(D^Lb(s,\cdot,\cdot)(\sL_{X_t^\mu})(X_t^\mu))(\ti z+\th(\ti y-\ti z))|^2\d\th\bigg)^\ff 1 2\big|_{y=X_t^{\mu_{\ve,\phi}},z=X_s^{\mu_{\ve,\phi}},\ti y=X_t^\mu,\ti z=X_s^\mu}\cr
&\qquad\times|X_t^{\mu_{\ve,\phi}}-X_s^{\mu_{\ve,\phi}}|,\cr
\Lambda_3&:=(t-s)^H\bigg(\E\int_0^1\bigg(\E|D^L(D^Lb(s,z,\cdot)(\cdot)(u))(\sL_{X^{\ve,\phi}_{s,t}(\th)})(X^{\ve,\phi}_{s,t}(\th))\cr
&\qquad\qquad-D^L(D^Lb(s,\ti z,\cdot)(\cdot)(v))(\sL_{X_{s,t}(\th)})(X_{s,t}(\th))|^2\bigg)\big|_{u=X_t^{\mu_{\ve,\phi}},v=X_t^\mu}
\d\th\bigg)^\ff 1 2\big|_{z=X_s^{\mu_{\ve,\phi}},\ti z=X_s^\mu},\cr
\Lambda_4&:=\bigg(\E\bigg(\int_0^1|\na(D^Lb(s,z,\cdot)(\sL_{X_s^{\mu_{\ve,\phi}}})(\cdot))(X^{\ve,\phi}_{s,t}(\th))\cr
&\qquad\qquad-\na(D^Lb(s,\ti z,\cdot)(\sL_{X_s^\mu})(\cdot))(X_{s,t}(\th))|^2\d\th\cdot|X_t^{\mu_{\ve,\phi}}-X_s^{\mu_{\ve,\phi}}|^2\bigg)\bigg)^\ff 1 2|_{z=X_s^{\mu_{\ve,\phi}},\ti z=X_s^\mu},
\end{align*}
and recall that for any $\th\in[0,1]$, $X_{s,t}(\th)=X_s^\mu+\th(X_t^\mu-X_s^\mu)$.\\
Note that due to \eqref{Pf2(ND)Le}, it follows that as $\ve$ or $\|\phi\|_{L^2_\mu}$ goes to zero,
$X_s^{\mu_{\ve,\phi}}$ and $X^{\ve,\phi}_{s,t}(\th)$ converge respectively to $X_s^\mu$ and $X_{s,t}(\th)$ in probability for any $s,t\in[0,T]$ and $\th\in[0,1]$.
Then, using (A3) again and applying the dominated convergence theorem, we deduce that
\beg{align}\label{Pf16-(ND)Th}
\lim_{\ve\ra0}&\E\left(\int_0^T\E\left(t^{2H-1}\left(\int_0^t\ff {|\Theta_6(t,s)|s^{\ff 1 2-H}}{(t-s)^{\ff 1 2+H}}\d s\right)^2\big|\sF_0\right)\d t\right)^\ff 1 2=0
\end{align}
and
\beg{align}\label{Pf21-(ND)Th}
\lim_{\|\phi\|_{L^2_\mu}\ra 0}\sup\limits_{\ve\in(0,1]}\ff{\E\left(\int_0^T\E\left(t^{2H-1}\left(\int_0^t\ff {|\Theta_6(t,s)|s^{\ff 1 2-H}}{(t-s)^{\ff 1 2+H}}\d s\right)^2\big|\sF_0\right)\d t\right)^\ff 1 2}{\|\phi\|_{L^2_\mu}}=0.
\end{align}
Hence, combining these and \eqref{Pf11-(ND)Th}-\eqref{Pf14-(ND)Th}, \eqref{Pf15-(ND)Th} with \eqref{Pf17-(ND)Th}, and applying Remark \ref{(ND)Re1}, we conclude that
\beg{align*}
\lim_{\ve\ra0}\E\left(\int_0^T\E(|J_4(t)|^2|\sF_0)\d t\right)^\ff 1 2=0
\end{align*}
and
\beg{align*}
\lim_{\|\phi\|_{L^2_\mu}\ra 0}\sup\limits_{\ve\in(0,1]}\ff{\E\left(\int_0^T\E(|J_4(t)|^2|\sF_0)\d t\right)^\ff 1 2}{\|\phi\|_{L^2_\mu}}=0.
\end{align*}
In conjunction with  \eqref{Pf18-(ND)Th}, \eqref{ForRem}, \eqref{Pf8-(ND)Th} and \eqref{Pf9-(ND)Th}, the above two inequalities imply
\beg{align*}
\lim_{\ve\ra 0^+}\E|\delta(h^{\ve,\phi})-\delta(h^\phi)|=0,\ \ \
\lim_{\|\phi\|_{L^2_\mu}\ra 0}\sup\limits_{\ve\in(0,1]}\ff{\E|\delta(h^{\ve,\phi})-\delta(h^\phi)|}{\|\phi\|_{L^2_\mu}}=0.
\end{align*}

Therefore, the assertions follow from Theorem \ref{(Ge)Th1}.
\qed

We conclude this part with a remark.
\beg{rem}\label{Non-Re2}
(i) Compared with the relevant result on distribution dependent SDE driven by the standard Brownian motion ($H=\ff 1 2$) shown in \cite[Theorem 2.1]{RW},
it is easy to see that our above result Theorem \ref{(ND)Th} applies to more general SDEs since we replace $B^\ff 1 2$ with fractional Brownian motion $B^H$ with arbitrary $H\in(\ff 1 2,1)$ as driving process.
Furthermore, due to the appearance of $J_4(t)$ in \eqref{ForRem}, essential difficulties are overcome in the analysis of Bismut formula for the $L$-derivative.

(ii) Combining the above proof and Remark \ref{Non-Re1}, we can derive the estimate of the $L$-derivative as the following:
\beg{align}\label{Non-Re2-1}
&\|D^L(P_T f)(\mu)\|=\sup_{\|\phi\|_{L^2_\mu}\le1}|D^L_\phi(P_T f)(\mu)|\cr
&\le\left(a_1(T,K,\tilde{K},H)+\ff {a_2(T,K,\tilde{K},H)} {T^H}\right)\left[(P_Tf^2)(\mu)-(P_T f(\mu))^2\right]^{\ff 1 2},\ \ f\in\sB(\R^d),
\end{align}
where $a_i(T,K,\tilde{K},H), i=1,2$ are two positive constants satisfying
\beg{align*}
a_1(T,K,\tilde{K},H)+\ff {a_2(T,K,\tilde{K},H)} {T^H}=O\left(\ff 1 {T^H}\right)\ \mathrm{when}\ T\ra0.
\end{align*}
Indeed, according to Theorem \ref{(ND)Th} and the H\"{o}lder inequality, we have
\beg{align*}
|D^L_\phi(P_T f)(\mu)|^2&=\left[\E\left(f(X_T^\mu)\int_0^T\<K_H^{-1}(R_H h^\phi)(t),\d W_t\>\right)\right]^2\cr
&=\left[\E\left((f(X_T^\mu)-P_T f(\mu))\int_0^T\<K_H^{-1}(R_H h^\phi)(t),\d W_t\>\right)\right]^2\cr
&\le\left[(P_Tf^2)(\mu)-(P_T f(\mu))^2\right]\int_0^T\E|K_H^{-1}(R_H h^\phi)(t)|^2\d t.
\end{align*}
Then, along the same lines as in \eqref{Pf2-(ND)Th}, \eqref{Pf4-(ND)Th}, \eqref{Pf5'(ND)Th} and \eqref{Add1-I5},
applying Remark \ref{Non-Re1} and taking into account of the relation $\sup_{s\in[0,T]}\E|\varrho(s)|^2\leq C(\ff 1 T+1)^2\|\phi\|_{L^2_\mu}^2$,
we obtain the estimate \eqref{Non-Re2-1}.\\
In addition, following the same arguments as in the proof of \cite[Corollary 2.2 (2)]{RW} and using \eqref{Non-Re2-1},
we give the total variation distance estimate for the difference between $\sL_{X_T^\mu}$ and $\sL_{X_T^\nu}$ with different initial distributions $\mu$ and $\nu$:
\beg{align*}
\|\sL_{X_T^\mu}-\sL_{X_T^\nu}\|_{\mathrm{var}}:=\sup_{A\in\sB(\R^d)}|\sL_{X_T^\mu}(A)-\sL_{X_T^\nu}(A)|\le C(T,K,\tilde{K},H)\W_2(\mu,\nu),\ \ \mu,\nu\in\sP_2(\R^d).
\end{align*}
\end{rem}

\subsection{Bismut formula: the degenerate case}

Let $d=m+l, b=(b^{(1)},b^{(2)})$ and $B^H$ be a $l$-dimensional fractional Brownian motion.
We now consider the following distribution dependent degenerate SDE:
\begin{equation}\label{(Deg)eq-1}
\begin{cases}
 \textnormal\d X^{(1)}_t=b^{(1)}(t,X_t)\d t,\\
 \textnormal\d X^{(2)}_t=b^{(2)}(t,X_t,\sL_{X_t})\d t+\sigma(t)\d B_t^H,
\end{cases}
\end{equation}
where $X_t=(X^{(1)}_t,X^{(2)}_t), b^{(1)}:[0,T]\times\R^{m+l}\ra\R^m, b^{(2)}:[0,T]\times\R^{m+l}\times\sP_2(\R^{m+l})\ra\R^l, \si(t)$ is an invertible $l\times l$-matrix for every $t\in[0,T]$.
It is obvious that \eqref{(Deg)eq-1} can be rewritten as follows
\beg{align}\label{+(Deg)eq-1}
\d X_t=(b^{(1)}(t,X_t),b^{(2)}(t,X_t,\sL_{X_t}))\d t+(0,\sigma(t)\d B_t^H).
\end{align}
Let us mention that, as in the Brownian motion case (see, e.g., \cite{BRW,RW}), when taking the special choices of $b^{(1)}, b^{(2)}$ and $\si$,
the above model will reduce to distribution dependent stochastic Hamiltonian system with fractional noise.

In the current part, we aim to establish the Bismut formula for the $L$-derivative of \eqref{+(Deg)eq-1} with the help of Theorem \ref{(Ge)Th1}.
To this end, we will impose the following condition.
\beg{enumerate}
\item[(C1)]For every $t\in[0,T], b^{(1)}(t,\cdot)\in C^1(\R^{m+l}\ra\R^m), b^{(2)}(t,\cdot,\cdot)\in C^{1,(1,0)}(\R^{m+l}\times\sP_2(\R^{m+l})\ra\R^l)$.
Moreover, there exists a constant $K>0$ such that for any $ t\in[0,T],\ x,y\in\R^{m+l},\ \mu\in\sP_2(\R^{m+l})$,
\beg{align*}
\|\na b^{(1)}(t,\cdot)(x)\|+\|\na b^{(2)}(t,\cdot,\mu)(x)\|+|D^Lb^{(2)}(t,x,\cdot)(\mu)(y)|\le K,
\end{align*}
and $\sup_{t\in[0,T]}(|b^{(1)}(t,0)|+|b^{(2)}(t,0,\de_0)|+\|\si(t)\|)\leq K$.
\end{enumerate}
It is easily checked that (C1) implies (H), so that there exists a unique solution to \eqref{+(Deg)eq-1} with any initial value $X_0\in L^2(\Omega\ra\R^{m+l},\sF_0,\P)$ thanks to Theorem \ref{(EX)The1}.
For any $\ve\in[0,1]$ and  $\phi\in L^2(\R^{m+l}\ra\R^{m+l},\mu)$, denote $X_t^{\mu_{\varepsilon,\phi}}$ by the solution of \eqref{+(Deg)eq-1} with $X_0^{\mu_{\varepsilon,\phi}}=(\mathrm{Id}+\varepsilon\phi)(X_0)$ (as before, in order to ease notation, we simply write $\mu_{\varepsilon,\phi}=\sL_{(\mathrm{Id}+\varepsilon\phi)(X_0)}$).
For any $s_0\in[0,T)$, let $\{g(t)\}_{t\in[0,T]}=\{(g^{(1)}(t),g^{(2)}(t))\}_{t\in[0,T]}$ be a stochastic process on $\R^{m+l}$ with differentiable paths satisfying
\beg{align}\label{(Deg)for-1}
g^{(1)}(t)&=\na_{\phi(X_0)} X_{t\wedge s_0}^{\mu_{\ve,\phi},(1)}+\int_{t\wedge s_0}^t\na_{g(s)}b^{(1)}(s,\cdot)(X_s^{\mu_{\ve,\phi}})\d s,\ \ t\in[0,T],\\
g^{(2)}(t)&=\na_{\phi(X_0)} X_{t}^{\mu_{\ve,\phi},(2)}, \ \ t\in[0,s_0],\label{(Deg)for-3}
\end{align}
and then put for each $t\in[0,T]$,
\beg{align}\label{(Deg)for-2}
(R_H h_{s_0}^{\ve,\phi})(t)
&:=\int_{t\wedge s_0}^t\si^{-1}(s)\left[\na_{g(s)}b^{(2)}(s,\cdot,\sL_{X_s^{\mu_{\varepsilon,\phi}}})(X_s^{\mu_{\varepsilon,\phi}})\right.\nonumber\\
&\quad\left.+\left(\E\<D^Lb^{(2)}(s,y,\cdot)(\sL_{X_s^{\mu_{\varepsilon,\phi}}})(X_s^{\mu_{\varepsilon,\phi}}), \na_{\phi(X_0)} X_s^{\mu_{\varepsilon,\phi}}\>\right)|_{y=X_s^{\mu_{\varepsilon,\phi}}}-(g^{(2)})'(s)\right]\d s.
\end{align}

\beg{prp}\label{(Deg)Th1}
Assume that (C1) holds, and that for any  $\ve\in[0,1]$ and $s_0\in[0,T)$,
let $g=(g^{(1)},g^{(2)})$ and $h_{s_0}^{\ve,\phi}\in\mathrm{Dom}\delta\cap\cH$ be respectively given in \eqref{(Deg)for-1}, \eqref{(Deg)for-3} and \eqref{(Deg)for-2} such that $g(T)=0, \int_0^1\left(\E\delta^2(h_{s_0}^{\tau,\phi})\right)^\ff 1 2\d\tau<\infty$ and
$$\lim_{\ve\ra 0^+}\E|\delta(h^{\ve,\phi})-\delta(h^\phi)|=0,\ \  \phi\in L^2(\R^{m+l}\ra\R^{m+l},\mu).$$
Then we have\\
(i) For any $f\in\sB_b(\R^{m+l})$, $P_Tf$ is weakly $L$-differentiable at $\mu$, and moreover
\beg{align*}
D^L_\phi(P_T f)(\mu)=\E(f(X_T^\mu)\delta(h^\phi)),\ \  \phi\in L^2(\R^{m+l}\ra\R^{m+l},\mu).
\end{align*}
(ii) For any $f\in\sB_b(\R^{m+l}), P_T f$ is $L$-differentiable at $\mu$, if $\E\delta^2(h^\phi)\leq \tilde{L}\|\phi\|^2_{L^2_\mu}$ with a constant $\tilde{L}>0$ and
\beg{align*}
 \lim_{\|\phi\|_{L^2_\mu}\ra 0}\sup\limits_{\ve\in(0,1]}\ff{\E|\delta(h^{\ve,\phi})-\delta(h^\phi)|}{\|\phi\|_{L^2_\mu}}=0.
\end{align*}
\end{prp}

\beg{proof}
Owing to \eqref{(Deg)for-2}, it is easy to see that $(R_H h_{s_0}^{\ve,\phi})(t)=0$ for all $t\in[0,s_0]$.
Now, to show assertions (i) and (ii), by Theorem \ref{(Ge)Th1} it remains to verify
$\D_{R_H h_{s_0}^{\ve,\phi}}X_T^{\mu_{\varepsilon,\phi}}=\na_{\phi(X_0)} X_T^{\mu_{\varepsilon,\phi}}$.

On one hand, according to Proposition \ref{Prop-Mall}, we obtain that for any $\ve\in[0,1],\phi\in L^2(\R^{m+l}\ra\R^{m+l},\mu)$ and $s_0\in[0,T)$,
the Malliavin derivative process $(Y_t:=\D_{R_H h_{s_0}^{\ve,\phi}}X_t^{\mu_{\ve,\phi}})_{t\in[0,T]}$ solves
\beg{align}\label{(Deg)Mall-1}
Y_t&=\int_0^t\left(\na_{Y_s}b^{(1)}(s,\cdot)(X_s^{\mu_{\varepsilon,\phi}}),
\na_{Y_s}b^{(2)}(s,\cdot,\sL_{X_s^{\mu_{\varepsilon,\phi}}})(X_s^{\mu_{\varepsilon,\phi}})\right)\d s\nonumber\\
&\quad+\int_0^t\left(0,\si(s)\d(R_H h_{s_0}^{\ve,\phi})(s)\right)\nonumber\\
&=\int_{t\wedge s_0}^t\left(\na_{Y_s}b^{(1)}(s,\cdot)(X_s^{\mu_{\varepsilon,\phi}}),
\na_{Y_s}b^{(2)}(s,\cdot,\sL_{X_s^{\mu_{\varepsilon,\phi}}})(X_s^{\mu_{\varepsilon,\phi}})\right)\d s\nonumber\\
&\quad+\int_{t\wedge s_0}^t\left(0,\si(s)\d(R_H h_{s_0}^{\ve,\phi})(s)\right),
\end{align}
where the second equality is due to the fact that $(R_H h_{s_0}^{\ve,\phi})(t)=0$ for $t\in[0,s_0]$.\\
Observe that by \eqref{(Deg)for-3} and \eqref{(Deg)for-2}, we have that for any $t\in[0,T]$,
\beg{align*}
g^{(2)}(t)
&=\na_{\phi(X_0)} X_{t\wedge s_0}^{\mu_{\varepsilon,\phi},(2)}+\int_{t\wedge s_0}^t\left[\na_{g(s)}b^{(2)}(s,\cdot,\sL_{X_s^{\mu_{\varepsilon,\phi}}})(X_s^{\mu_{\varepsilon,\phi}})\right.\cr
&\qquad\qquad\qquad\qquad+\left.\left(\E\<D^Lb^{(2)}(s,y,\cdot)(\sL_{X_s^{\mu_{\varepsilon,\phi}}})(X_s^{\mu_{\varepsilon,\phi}}), \na_{\phi(X_0)} X_s^{\mu_{\varepsilon,\phi}}\>\right)|_{y=X_s^{\mu_{\varepsilon,\phi}}}\right]\d s\cr
&\quad-\int_{t\wedge s_0}^t\si(s)\d(R_H h_{s_0}^{\ve,\phi})(s).
\end{align*}
Then, combining this with \eqref{(Deg)for-1} implies that for any $t\in[0,T]$,
\beg{align*}
g(t)&=\na_{\phi(X_0)} X_{t\wedge s_0}^{\mu_{\varepsilon,\phi}}
+\int_{t\wedge s_0}^t\left[\left(\na_{g(s)}b^{(1)}(s,\cdot)(X_s^{\mu_{\ve,\phi}}),\na_{g(s)}b^{(2)}(s,\cdot,\sL_{X_s^{\mu_{\varepsilon,\phi}}})(X_s^{\mu_{\varepsilon,\phi}})\right)\right.\cr
&\qquad\qquad\qquad\quad+\left.\left(0,\left(\E\<D^Lb^{(2)}(s,y,\cdot)(\sL_{X_s^{\mu_{\varepsilon,\phi}}})(X_s^{\mu_{\varepsilon,\phi}}), \na_{\phi(X_0)} X_s^{\mu_{\varepsilon,\phi}}\>\right)|_{y=X_s^{\mu_{\varepsilon,\phi}}}\right)\right]\d s\cr
&\quad-\int_{t\wedge s_0}^t\left(0,\si(s)\d(R_H h_{s_0}^{\ve,\phi})(s)\right).
\end{align*}
Consequently, in conjunction with \eqref{(Deg)Mall-1}, the above relation yields that for any $t\in[0,T]$,
\beg{align}\label{Pf1(De)Th1}
Y_t+g(t)&=\na_{\phi(X_0)} X_{t\wedge s_0}^{\mu_{\varepsilon,\phi}}
+\int_{t\wedge s_0}^t\left[\left(\na_{Y_s+g(s)}b^{(1)}(s,\cdot)(X_s^{\mu_{\ve,\phi}}),\na_{Y_s+g(s)}b^{(2)}(s,\cdot,\sL_{X_s^{\mu_{\varepsilon,\phi}}})(X_s^{\mu_{\varepsilon,\phi}})\right)\right.\nonumber\\
&\qquad\qquad\qquad+\left.\left(0,\left(\E\<D^Lb^{(2)}(s,y,\cdot)(\sL_{X_s^{\mu_{\varepsilon,\phi}}})(X_s^{\mu_{\varepsilon,\phi}}), \na_{\phi(X_0)} X_s^{\mu_{\varepsilon,\phi}}\>\right)|_{y=X_s^{\mu_{\varepsilon,\phi}}}\right)\right]\d s.
\end{align}

On the other hand, for any $\ve\in[0,1],\phi\in L^2(\R^{m+l}\ra\R^{m+l},\mu)$ and $s_0\in[0,T)$,
applying Proposition \ref{Prop-ParV} with $\eta=\phi(X_0)$,
one sees that the directional derivative process $(Z_t:=\na_{\phi(X_0)} X_t^{\mu_{\ve,\phi}})_{t\in[0,T]}$ solves
\beg{align}\label{(Deg)ParV-1}
Z_t&=\phi(X_0)+\int_0^t\left[\left(\na_{Z_s}b^{(1)}(s,\cdot)(X_s^{\mu_{\ve,\phi}}),\na_{Z_s}b^{(2)}(s,\cdot,\sL_{X_s^{\mu_{\ve,\phi}}})(X_s^{\mu_{\ve,\phi}})\right)\right.\nonumber\\
&\qquad\qquad\qquad\quad\left.+\left(0,\left(\E\<D^Lb^{(2)}(s,y,\cdot)(\sL_{X_s^{\mu_{\ve,\phi}}})(X_s^{\mu_{\ve,\phi}}), Z_s\>\right)|_{y=X_s^{\mu_{\ve,\phi}}}\right)\right]\d s\nonumber\\
&=Z_{t\wedge s_0}+\int_{t\wedge s_0}^t\left[\left(\na_{Z_s}b^{(1)}(s,\cdot)(X_s^{\mu_{\ve,\phi}}),\na_{Z_s}b^{(2)}(s,\cdot,\sL_{X_s^{\mu_{\ve,\phi}}})(X_s^{\mu_{\ve,\phi}})\right)\right.\nonumber\\
&\qquad\qquad\qquad\quad\left.+\left(0,\left(\E\<D^Lb^{(2)}(s,y,\cdot)(\sL_{X_s^{\mu_{\ve,\phi}}})(X_s^{\mu_{\ve,\phi}}), Z_s\>\right)|_{y=X_s^{\mu_{\ve,\phi}}}\right)\right]\d s.
\end{align}
Owing to \eqref{Pf1(De)Th1}, \eqref{(Deg)ParV-1} and the uniqueness of solutions of the ODE, we conclude that
\beg{align*}
\D_{R_H h_{s_0}^{\ve,\phi}}X_t^{\mu_{\varepsilon,\phi}}+g(t)=\na_{\phi(X_0)} X_t^{\mu_{\varepsilon,\phi}}, \ \  t\in[0,T].
\end{align*}
This implies $\D_{R_H h_{s_0}^{\ve,\phi}}X_T^{\mu_{\varepsilon,\phi}}=\na_{\phi(X_0)} X_T^{\mu_{\varepsilon,\phi}}$ due to $g(T)=0$.
The proof is now complete.
\end{proof}
Next, we intend to apply Proposition \ref{(Deg)Th1} with concrete choices of $g=(g^{(1)},g^{(2)})$.
Without of lost generality, we consider a special case of \eqref{(Deg)eq-1} with $b^{(1)}(t,x,y)=Ax+By$,
where $A$ and $B$ are two matrices of order $m\times m$ and $m\times l$, respectively.
That is,
\begin{equation}\label{(Deg)eq-2}
\begin{cases}
 \textnormal\d X^{(1)}_t=(AX^{(1)}_t+B X^{(2)}_t)\d t,\\
 \textnormal\d X^{(2)}_t=b^{(2)}(t,X_t,\sL_{X_t})\d t+\sigma(t)\d B_t^H.
\end{cases}
\end{equation}
For the equation \eqref{(Deg)eq-2}, we impose additional conditions on $b^{(2)}$ and $\si$ which are similar to (A2) and (A3).
\beg{enumerate}
\item[(C2)] There exists a constant $\tilde{K}>0$ such that

\item[(i)]  for any $t,s\in[0,T],\ x,y,z_1,z_2\in\R^{m+l},\ \mu,\nu\in\sP_2(\R^{m+l})$,
\beg{align*}
&\|\na b^{(2)}(t,\cdot,\mu)(x)-\na b^{(2)}(s,\cdot,\nu)(y)\|+|D^Lb^{(2)}(t,x,\cdot)(\mu)(z_1)-D^Lb^{(2)}(s,y,\cdot)(\nu)(z_2)|\cr
&\le \tilde{K}(|t-s|^{\alpha_0}+|x-y|^{\beta_0}+|z_1-z_2|^{\gamma_0}+\W_2(\mu,\nu)),
\end{align*}
where $\alpha_0\in(H-1/2,1]$ and $\beta_0,\gamma_0\in(1-1/(2H),1]$.

\item[(ii)]  $\si^{-1}$ is H\"{o}lder continuous of order $\delta_0\in(H-1/2,1]$:
\beg{align*}
\|\si^{-1}(t)-\si^{-1}(s)\|\le\tilde{K}|t-s|^{\delta_0}, \ \ t,s\in[0,T].
\end{align*}

\item[(C3)] The derivatives
$$\partial_t(D^Lb^{(2)}(\cdot,x,\cdot)(\mu)(y))(t),\ \na(D^Lb^{(2)}(t,\cdot,\cdot)(\mu)(y))(x),$$
$$D^L(D^Lb^{(2)}(t,x,\cdot)(\cdot)(y))(\mu)(z),\ \na(D^Lb^{(2)}(t,x,\cdot)(\mu)(\cdot))(y)$$
exist and are bounded continuous in the corresponding arguments $(t,x,\mu,y)$ or $(t,x,\mu,y,z)$.
We denote the bounded constants by a common one $\bar K>0$.
\end{enumerate}
For any $s_0\in[0,T)$, let
\beg{align}\label{Pf3(De)Th2}
U^{s_0}_t=\int_{s_0}^t\ff{(s-s_0)(T-s)}{T^2} \e^{(T-s)A}BB^*\e^{(T-s)A^*}\d s\geq\rho(t)\mathrm{I}_{m\times m}, \ t\in(s_0,T],
\end{align}
where $\rho\in C([0,T])$ satisfies $\rho(t)>0$ for any $t\in(0,T]$ and $\mathrm{I}_{m\times m}$ is the identity matrix on $\R^m\times\R^m$, and set for $t\in[0,T]$,
\beg{align}\label{Pf1(De)Th2}
g^{(1)}(t):=\e^{(t-t\wedge s_0)A}\na_{\phi(X_0)} X_{t\wedge s_0}^{\mu_{\ve,\phi},(1)}+\int_{t\wedge s_0}^t\e^{(t-s)A}Bg^{(2)}(s)\d s
\end{align}
and
\beg{align}\label{Pf2(De)Th2}
g^{(2)}(t)&:=\ff {T-t} {T-t\wedge s_0}\na_{\phi(X_0)} X_{t\wedge s_0}^{\mu_{\ve,\phi},(2)}\nonumber\\
&\quad-\ff{(t-t\wedge s_0)(T-t)}{T^2}B^*\e^{(T-t)A^*}(U^{s_0}_T)^{-1}\e^{(T-s_0)A}\na_{\phi(X_0)} X_{t\wedge s_0}^{\mu_{\ve,\phi},(1)}\nonumber\\
&\quad-\ff{(t-t\wedge s_0)(T-t)}{T^2}B^*\e^{(T-t)A^*}(U^{s_0}_T)^{-1}\int_{s_0}^T\ff {T-s} {T-s_0}\e^{(T-s)A}B\na_{\phi(X_0)} X_{s_0}^{\mu_{\ve,\phi},(2)}\d s.
\end{align}
Then our main result in the current part can be stated in the following theorem.

\beg{thm}\label{(Deg)Th2}
Assume that (C1), (C2) and (C3) hold.
Then for any $\mu\in\sP_2(\R^{m+l})$ and $f\in\sB_b(\R^{m+l}), P_T f$ is $L$-differentiable at $\mu$ such that
\beg{align*}
D^L_\phi(P_T f)(\mu)=\E\left(f(X_T^\mu)\int_0^T\<K_H^{-1}(R_H h^\phi)(t),\d W_t\>\right),\ \ \phi\in L^2(\R^{m+l}\ra\R^{m+l},\mu),
\end{align*}
where $h^\phi\in\mathrm{Dom}\delta\cap\cH$ and satisfies for every $t\in[0,T]$,
\beg{align*}
(R_H h^{\phi})(t)
&=\int_0^t\si^{-1}(s)\left[\na_{g(s)}b^{(2)}(s,\cdot,\sL_{X_s^{\mu}})(X_s^{\mu})\right.\cr
&\qquad\qquad\qquad\left.+\left(\E\<D^Lb^{(2)}(s,y,\cdot)(\sL_{X_s^{\mu}})(X_s^{\mu}), \na_{\phi(X_0)}
X_s^{\mu}\>\right)|_{y=X_s^{\mu}}-(g^{(2)})'(s)\right]\d s
\end{align*}
and $g=(g^{(1)},g^{(2)})$ is given by \eqref{Pf1(De)Th2} and \eqref{Pf2(De)Th2} for $s_0=0$ and
$(\na_{\phi(X_0)} X_\cdot^{\mu,(1)},\na_{\phi(X_0)} X_\cdot^{\mu,(2)})$
replacing $(\na_{\phi(X_0)} X_\cdot^{\mu_{\ve,\phi},(1)},\na_{\phi(X_0)} X_\cdot^{\mu_{\ve,\phi},(2)})$.
\end{thm}

\beg{proof}
Observe first that by \eqref{Pf3(De)Th2}, we can see that $U^{s_0}_t$ is invertible with $\|(U^{s_0}_t)^{-1}\|\leq 1/ \rho(t)$ for every $t\in(s_0,T]$,
and then $g^{(1)}(t)$ and $g^{(2)}(t)$ given respectively by \eqref{Pf1(De)Th2} and \eqref{Pf2(De)Th2} are well-defined.
Since $\na b^{(1)}=(A,B)$, we have that $g^{(1)}(t)$ satisfies \eqref{(Deg)for-1}.
Owing to \eqref{Pf2(De)Th2}, it is readily checked that $g^{(2)}(T)=0$ and \eqref{(Deg)for-3} holds.
Besides, by \eqref{Pf1(De)Th2} and \eqref{Pf2(De)Th2} again, we have
\beg{align*}
&g^{(1)}(T)=\e^{(T-s_0)A}\na_{\phi(X_0)} X_{s_0}^{\mu_{\ve,\phi},(1)}+\int_{s_0}^T\e^{(T-s)A}Bg^{(2)}(s)\d s\cr
&=\e^{(T-s_0)A}\na_{\phi(X_0)} X_{s_0}^{\mu_{\ve,\phi},(1)}+\int_{s_0}^T\ff {T-s} {T-s_0}\e^{(T-s)A}B\na_{\phi(X_0)} X_{s_0}^{\mu_{\ve,\phi},(2)}\d s\cr
&\quad-\int_{s_0}^T\ff{(s-s_0)(T-s)}{T^2}\e^{(T-s)A}BB^*\e^{(T-s)A^*}\d s(U^{s_0}_T)^{-1}\e^{(T-s_0)A}\na_{\phi(X_0)} X_{s_0}^{\mu_{\ve,\phi},(1)}\cr
&\quad-\int_{s_0}^T\ff{(s-s_0)(T-s)}{T^2}\e^{(T-s)A}BB^*\e^{(T-s)A^*}\d s(U^{s_0}_T)^{-1}\int_{s_0}^T\ff {T-s} {T-s_0}\e^{(T-s)A}B\na_{\phi(X_0)} X_{s_0}^{\mu_{\ve,\phi},(2)}\d s\cr
&=0,
\end{align*}
where the last equality is due to the definition of $U^{s_0}_T$.\\
For any  $\ve\in[0,1]$ and $s_0\in[0,T)$, let
\beg{align}\label{Pf4(De)Th2}
\tilde{h}_{s_0}^{\ve,\phi}(t)
&=\int_{t\wedge s_0}^t\si^{-1}(s)\left[\na_{g(s)}b^{(2)}(s,\cdot,\sL_{X_s^{\mu_{\varepsilon,\phi}}})(X_s^{\mu_{\varepsilon,\phi}})\right.\cr
&\left.+\left(\E\<D^Lb^{(2)}(s,y,\cdot)(\sL_{X_s^{\mu_{\varepsilon,\phi}}})(X_s^{\mu_{\varepsilon,\phi}}), \na_{\phi(X_0)}
X_s^{\mu_{\varepsilon,\phi}}\>\right)|_{y=X_s^{\mu_{\varepsilon,\phi}}}-(g^{(2)})'(s)\right]\d s,
\end{align}
where  $g^{(1)}$ and  $g^{(2)}$ are defined by  \eqref{Pf1(De)Th2} and \eqref{Pf2(De)Th2}, respectively.
Noting that the right-hand side of \eqref{Pf4(De)Th2} belongs to $I_{0+}^{H+\ff 1 2}(L^2([0,T],\R^l))$ due to (C1) and (C2),
there exists $h_{s_0}^{\ve,\phi}\in\cH$ such that $R_H h_{s_0}^{\ve,\phi}=\tilde{h}_{s_0}^{\ve,\phi}$.
Under (C1) and (C2), the adaptation of our calculations in Lemma \ref{(ND)Le} to the present degenerate case are straightforward.
Then resorting to the same techniques as in \eqref{addfor-de} and \cite[Theorem 3.2]{Fan19},
we obtain that  $K_H^{-1}(R_H h_{s_0}^{\ve,\phi})=K_H^*h_{s_0}^{\ve,\phi}\in L^2_a([0,T]\times\Omega,\R^l)$
(which implies that $\delta(h_{s_0}^{\ve,\phi})=\delta_W(K_H^*h_{s_0}^{\ve,\phi})=\int_0^T\<K_H^{-1}(R_H h_{s_0}^{\ve,\phi})(t),\d W_t\>$),
and
$\int_0^1\left(\E\delta^2(h_{s_0}^{\tau,\phi})\right)^\ff 1 2\d\tau<\infty$
as well as $\E\delta^2(h^\phi)\leq C_{T,K,\tilde{K},H}\|\phi\|^2_{L^2_\mu}$.
Furthermore, with the help of (C3) and by a similar analysis of \eqref{Pf18-(ND)Th} we derive that
\beg{align*}
\lim_{\ve\ra 0^+}\E|\delta(h^{\ve,\phi})-\delta(h^\phi)|=0 \ \ \mathrm{and}\ \ \lim_{\|\phi\|_{L^2_\mu}\ra 0}\sup\limits_{\ve\in(0,1]}\ff{\E|\delta(h^{\ve,\phi})-\delta(h^\phi)|}{\|\phi\|_{L^2_\mu}}=0.
\end{align*}

Therefore, the assertions follow from Proposition \ref{(Deg)Th1}.
\end{proof}

\beg{rem}\label{De-Re1}
Similar to Remark \ref{Non-Re2}(ii), we can obtain that there exist two positive constants $c_i(T,K,\tilde{K},H), i=1,2$ such that
\beg{align*}
&\|D^L(P_T f)(\mu)\|\le c(T,K,\tilde{K},H)\left[(P_Tf^2)(\mu)-(P_T f(\mu))^2\right]^{\ff 1 2}, \ \ f\in\sB(\R^{m+l}),\cr
&\|\sL_{X_T^\mu}-\sL_{X_T^\nu}\|_{\mathrm{var}}\le c(T,K,\tilde{K},H)\W_2(\mu,\nu),\ \ \mu,\nu\in\sP_2(\R^{m+l}),
\end{align*}
where
\beg{align*}
c(T,K,\tilde{K},H):=c_1(T,K,\tilde{K},H)+c_2(T,K,\tilde{K},H)\left(\ff {1} {T^H}+\ff 1{\rho(T)}+\ff {1}{T^H\rho(T)}\right)
\end{align*}
with $c(T,K,\tilde{K},H)=O\left(\ff 1 {T^H\rho(T)}\right)$ when $T\ra0$.
If the following Kalman rank condition
\beg{align*}
\mathrm{Rank}[B,AB,\cdots,A^kB]=m
\end{align*}
holds for some integer number $k\in[0,m-1]$, then \eqref{Pf3(De)Th2} is satisfied with $\rho(t)=\ff {C_1(t\wedge1)^{2(k+1)}}{T\e^{C_2T}}$
for two positive constants $C_i,i=1,2$ (see, e.g., \cite[Theorem 4.2]{WZ13}), which implies that
$c(T,K,\tilde{K},H)=O\left(\ff 1 {T^{2k+1+H}}\right)$ as $T\ra0$.
\end{rem}

\textbf{Acknowledgement}

X. Fan is supported by the DFG through the CRC 1283 Taming uncertainty and profiting from randomness and low
regularity in analysis, stochastics and their applications, and the Natural Science Foundation of Anhui Province (No. 2008085MA10).
X. Huang is supported by the National Natural Science Foundation of China (No. 11801406).


\begin{thebibliography}{17}
{\small

\setlength{\baselineskip}{0.14in}
\parskip=0pt

\bibitem{AKR96} S. Albeverio, Y. G. Kondratiev and M. R\"{o}ckner, Differential geometry of Poisson spaces,
\textit{C. R. Acad. Sci. Paris S\'{e}r. I Math.} {\bf 323} (1996), 1129--1134.

\bibitem{Alos&Mazet&Nualart01a} E. Al\`{o}s, O. Mazet and D. Nualart, Stochastic calculus with respect to Gaussian processes,
\textit{Ann. Probab.} {\bf 29} (2001), 766--801.


\bibitem{AN03} E. Al\`{o}s and D. Nualart, Stochastic integration with respect to the fractional Brownian motion,
\textit{Stochastics and Stochastic Reports} {\bf 75} (2003), 129--152.




\bibitem{Banos18} D. Ba\~{n}os, The Bismut-Elworthy-Li formula for mean-field stochastic differential equations, \textit{Ann. Inst. H. Poincar\'{e} Probab. Statist.} {\bf 54} (2018), 220--233.

\bibitem{BRW} J. Bao, P. Ren and F.-Y. Wang, Bismut formula for Lions derivative of distribution-path dependent SDEs,
\textit{J. Differential Equations} {\bf 282} (2021), 285--329.

\bibitem{Bismut84} J. Bismut, Large Deviation and The Malliavin Calculus, Birkh\"{a}user, Boston, MA, 1984.



\bibitem{BJ17} R. Buckdahn and S. Jing,  Mean-field SDE driven by a fractional Brownian motion and related stochastic control problem,
\textit{SIAM J. Control Optim.} {\bf 55} (2017), 1500--1533.

\bibitem{Buckdahn&L&Peng&ainer17a} R. Buckdahn, J. Li, S. Peng and C. Rainer,  Mean-field stochastic differential equations and
associated PDEs, \textit{Ann. Probab.} {\bf 2} (2017), 824--878.



\bibitem{Cardaliaguet13} P. Cardaliaguet, Notes on mean field games, P.-L. Lions lectures at Coll\`{e}ge de France,
https://www.ceremade.dauphine.fr/cardaliaguet/MFG20130420.pdf, 2013.


\bibitem{CD13}  R. Carmona and F. Delarue, Probabilistic analysis of mean-field games,
\textit{SIAM J. Control Optim.} {\bf 51} (2013), 2705--2734.

\bibitem{CD15} R. Carmona and F. Delarue, Forward-backward stochastic differential equations and controlled McKean-Vlasov dynamics,
\textit{Ann. Probab.} {\bf 43} (2015), 2647--2700.

\bibitem{Coh} D. L. Cohn, Measure Theory (Second Edition), Birkh\"{a}user, Boston, Mass., 2013.

\bibitem{Crisan&McMurray18a} D. Crisan and E. McMurray, Smoothing properties of McKean-Vlasov SDEs,
\textit{Probab. Theory Related Fields} {\bf 171} (2018), 97--148.



\bibitem{Decreusefond&Ustunel98a}L. Decreusefond and A. S. \"{U}st\"{u}nel, Stochastic analysis of the fractional Brownian motion,
\textit{Potential Anal.} {\bf 10} (1998), 177--214.


\bibitem{Elworthy&Li94} K. Elworthy and X. Li, Formulae for the derivatives of heat semigroups, \textit{J. Funct. Anal.} {\bf 125} (1994), 252--286.

\bibitem{Fan19} X. Fan, Derivative formulas and applications for degenerate stochastic differential equations with fractional noises,
\textit{J. Theoret. Probab.} {\bf 32} (2019), 1360--1381.

\bibitem{FR} X. Fan and Y. Ren, Bismut formulas and applications for stochastic (functional) differential equations
driven by fractional Brownian motions, \textit{Stoch. Dyn.} {\bf 17} (2017), 1750028.

\bibitem{HMC13}  M. Huang, R. Malham\'{e} and P. Caines, Large population stochastic dynamic games: closed-loop
McKean-Vlasov systems and the Nash certainty equivalence principle,
\textit{Commun. Inf. Syst.} {\bf 6} (2006), 221--251.

\bibitem{HRW19} X. Huang, M. R\"{o}ckner and F.-Y. Wang, Nonlinear Fokker-Planck equations for probability measures on path space and path-distribution dependent SDEs, \textit{Discrete Contin. Dyn. Syst.} {\bf 39} (2019), 3017--3035.


\bibitem{HW} X. Huang and F.-Y. Wang, Distribution dependent SDEs with singular coefficients, \textit{Stochastic Process.Appl.} {\bf 129} (2019), 4747--4770.

\bibitem{HW21} X. Huang and F.-Y.  Wang, Derivative estimates on distributions of McKean-Vlasov SDEs, \textit{Electron. J. Probab.} {\bf 26} (2021), 1--12.



\bibitem{Kac}  M. Kac, Foundations of kinetic theory,  in: Proceedings of the Third Berkeley Symposium on Mathematical
Statistics and Probability, University California Press, Berkeley, 1956, 171--197.



\bibitem{Li18a} J. Li,  Mean-field forward and backward SDEs with jumps and associated nonlocal quasi-linear
integral-PDEs, \textit{Stochastic Process. Appl.} {\bf 128} (2018), 3118--3180.

\bibitem{LWbook} W. Liu and M. R\"{o}ckner, Stochastic Partial Differential Equations: An Introduction, Universitext, Springer, 2015.


\bibitem{McKean}   H. P. McKean, Propagation of chaos for a class of non-linear parabolic equations,
in: Lecture Series in Differential Equations, Session 7, 1967, 41--57.


\bibitem{JMV01} J. M\'{e}min, Y. Mishura and E. Valkeila, Inequalities for the moments of Wiener integrals with respect to a fractional Brownian motion,
\textit{Statist. Probab. Lett.} {\bf 51} (2001), 197--206.




\bibitem{Nikiforov&Uvarov88} A. F. Nikiforov and V. B. Uvarov, Special Functions of Mathematical Physics, Birkh\"{a}user, Boston, 1988.

\bibitem{ND} D. Nualart, The Malliavin Calculus and Related Topics, Second edition, Springer-Verlag, Berlin, 2006.


\bibitem{RW} P. Ren and F.-Y. Wang, Bismut formula for Lions derivative of distribution dependent SDEs and applications,
\textit{J. Differential Equations} {\bf 267} (2019), 4745--4777.

\bibitem{RW21} P. Ren and F.-Y. Wang, Space-distribution PDEs for path independent additive functionals of McKean-Vlasov SDEs,
to appear in \textit{Infin. Dimens. Anal. Quantum Probab. Relat. Top.}



\bibitem{SK} S. G. Samko, A. A. Kilbas and O. I. Marichev, Fractional Integrals and Derivatives, Theory and Applications, Gordon and Breach Science Publishers, Yvendon, 1993.

\bibitem{Song} Y. Song,  Gradient estimates and exponential ergodicity for Mean-Field SDEs with jumps,
 \textit{J. Theoret. Probab.} {\bf 33} (2020), 201--238.

\bibitem{Stein70} E. M. Stein, Singular Integrals and Differentiability Properties of Functions, Princeton University Press, Princeton, 1970.

\bibitem{Sznitman}  A. Sznitman, Topics in propagation of chaos, in: Ecole d'Et\'{e} de Probabilit\'{e}s de Saint-Flour XIX-1989, Springer, Berlin, Heidelberg, 1991, 165--251.


\bibitem{Wang13a} F.-Y. Wang, Harnack Inequalities for Stochastic Partial Differential Equations, Springer, 2013.

\bibitem{Wang18}  F.-Y. Wang,  Distribution dependent SDEs for Landau type equations, \textit{Stochastic Process. Appl.} {\bf 128} (2018), 595--621.

\bibitem{WZ13}  F.-Y. Wang and X. Zhang, Derivative formula and applications for degenerate diffusion semigroups,
\textit{J. Math. Pures Appl.} {\bf 99} (2013), 726--740.






















}\end{thebibliography}
\end{document}